\documentclass[twoside]{article}
\usepackage{mor}
\usepackage{amsfonts}
\usepackage{amsmath}
\usepackage{amsxtra}
\received{}                              %
\revised{}                               %
\pubyear{0x}                             %
\pubmonth{Xxxxxxx}                       %
\volume{xx}                              %
\issue{x}                                %
\pages{xxx--xxx}                         %
\firstpage{0xxx}                         %
\DOI{10.1287/moor.xxxx.xxxx}             %
\startpagenumber{1}                      %

\title{A Limit Theorem for Financial Markets with Inert
Investors}

\ShortTitle{Investor Inertia}

\ShortAuthors{Bayraktar, Horst, Sircar}

\NumberOfAuthors{3}

\FirstAuthor{Erhan Bayraktar}

\FirstAuthorAddress{Department of Mathematics, University of
Michigan, 525 East University, Ann Arbor, MI 48109}
\FirstAuthorEmail{erhan@umich.edu} 

\SecondAuthor{Ulrich Horst} \SecondAuthorAddress{Department of
Mathematics, University of British Columbia,  1984 Mathematics
Road, Vancouver, BC, V6T 1Z2}
\SecondAuthorEmail{horst@math.ubc.ca}

\ThirdAuthor{Ronnie Sircar} \ThirdAuthorAddress{Department of
Operations Research \& Financial Engineering, Princeton
University, Princeton, NJ 08544}
\ThirdAuthorEmail{sircar@princeton.edu}


\MSCcodes{Primary:   60F13, 60G15,    ; Secondary:  91B28   } 

\ORMScodes{Primary: Probability   ; Secondary: Stochastic model applications  } 
\newcommand{\ma}{\mathbb{A}}
\newcommand{\n}{\mathbb{N}}
\newcommand{\re}{\mathbb{R}}
\newcommand{\E}{\mathbb{E}}
\newcommand{\Prob}{\mathbb{P}}

\newcommand{\mf}{\mathcal{F}}
\newcommand{\Var}{\textnormal{Var}}
\newcommand{\fbm}{fractional Brownian motion }
\renewcommand{\vec}[1]{\mbox{\boldmath $#1$}}
\newcommand{\eps}{\varepsilon}
\newcommand{\ind}{\textbf{1}}
\newenvironment{rmenumerate}
    {\begin{enumerate}}
    {\end{enumerate}}
\begin{document}
\maketitle

\begin{abstract}
We study the effect of investor inertia on stock price
fluctuations with a market microstructure model comprising many
small investors who are inactive most of the time. It turns out
that semi-Markov processes are tailor made for modelling inert
investors. With a suitable scaling, we show that when the price is
driven by the market imbalance, the log price process  is
approximated by a process with long range dependence and
non-Gaussian returns distributions, driven by a fractional
Brownian motion. Consequently, investor inertia may lead to
arbitrage opportunities for sophisticated market participants. The
mathematical contributions are a functional central limit theorem
for stationary semi-Markov processes, and approximation results
for stochastic integrals of continuous semimartingales with
respect to fractional Brownian motion.
\end{abstract}
\normalsize


\section{Introduction and Motivation}


We prove a functional central limit theorem for stationary
semi-Markov processes in which the limit process is a stochastic
integral with respect to fractional Brownian motion. Our
motivation is to develop a probabilistic framework within which to
analyze the aggregate effect of investor inertia on asset price
dynamics. We show that, in isolation, such infrequent trading
patterns can lead to long-range dependence in stock prices and
arbitrage opportunities for other more ``sophisticated'' traders.

\subsection{Market Microstructure Models for Financial Markets}

In mathematical finance, the dynamics of asset prices are usually
modelled by trajectories of some exogenously specified stochastic
process defined on some underlying probability space $(\Omega,
\mf, \Prob)$. Geometric Brownian motion has long become the
canonical reference model of financial price fluctuations. Since
prices are generated by the demand of market participants, it is
of interest to support such an approach by a microeconomic model
of interacting agents. 

In recent years there has been increasing interest in agent-based
models of financial markets. These models are capable of
explaining, often through simulations, many facts like the
emergence of herding behavior \cite{Lux}, volatility clustering
\cite{LuxMarchesi} or fat-tailed distributions of stock returns
\cite{CB} that are observed in financial data. Brock and Hommes
\cite{BrockHommes,BrockHommesb} proposed 
models with many traders where the asset price process is
described by {\em deterministic} dynamical systems. From numerical
simulations, they showed that financial price fluctuations can
exhibit chaotic behavior if the effects of technical trading
become too strong.

F\"ollmer and Schweizer \cite{Foellmer-Schweizer} took the
probabilistic point of view, with asset prices arising from a
sequence of temporary price equilibria in an exogenous random
environment of investor sentiment; see \cite{Foellmer94},
\cite{Horst02} or \cite{FHK} for similar approaches. Applying an
invariance principle to a sequence of suitably defined discrete
time models, they derived a diffusion approximation for the
logarithmic price process. Duffie and Protter
\cite{Duffie-Protter} also provided a mathematical framework for
approximating sequences of
stock prices by diffusion processes. 

All the aforementioned models assume that the agents trade the
asset in each period. 
At the end of each trading interval, the agents update their
expectations for the future evolution of the stock price and
formulate their excess demand for
the following period. 
However, small investors are not so efficient in their investment
decisions: they are typically inactive and actually trade only
occasionally. This may be because they are waiting to accumulate
sufficient capital to make further stock purchases; or they tend
to monitor their portfolios infrequently; or they are simply
scared of choosing the wrong investments; or they feel that as
long-term investors, they can defer action; or they put off the
time-consuming research necessary to make informed portfolio
choices. Long uninterrupted periods of inactivity may be viewed as
a form of investor inertia. The focus of this paper is the effect
of such investor inertia on asset prices in a model with
asynchronous order arrivals. See \cite{Kruk} for an alternative
micro-structure model with asynchronous trading.


\subsection{Inertia in Financial Markets}

Investor inertia is a common experience and is well documented.
The New York Stock Exchange (NYSE)'s survey of individual
shareownership in the United States, ``Shareownership2000''
\cite{NYSE}, demonstrates that many investors have very low levels
of trading activity. For example they find that ``23 percent of
stockholders with brokerage accounts report no trading at all,
while 35 percent report trading only once or twice in the last
year'' (see pages 58-59).  The NYSE survey (e.g. Table 28) also
reports that the average holding period for stocks is long, for
example 2.9 years in the early 90's.



Empirical evidence of inertia also appears in the economic
literature. For example, Madrian and Shea \cite{madrian} looked at
the reallocation of assets in employees' individual 401(k)
(retirement) plans\footnote{A 401k retirement plan is a special
type of account funded through pre-tax payroll deductions. The
funds in the account can be invested in a number of different
stocks, bonds, mutual funds or other assets, and are not taxed on
any capital gains, dividends, or interest until they are
withdrawn. The retirement savings vehicle was created by United
States Congress in 1981 and gets its name from the section of the
Internal Revenue Code that describes it.}
and found ``a status quo bias resulting from employee
procrastination in making or implementing an optimal savings
decision.'' A related study by Hewitt Associates (a management
consulting firm) found that in 2001, four out of five plan
participants did not do any trading in their 401(k)s. Madrian and
Shea explain that ``if the cost of gathering and evaluating the
information needed to make a 401(k) savings decision exceeds the
short-run benefit from doing so, individuals will procrastinate.''
The prediction of Prospect Theory \cite{kahneman} that investors
tend to hold onto losing stocks too long has also been observed
(\cite{shefrin}). 

A number of microeconomic models study investor caution with
regard to model risk, which is termed uncertainty aversion. Among
others, Dow and Werlang (\cite{Dow-Werlang}) and Simonsen and
Werlang (\cite{Simonsen-Werlang}) considered models of portfolio
optimization where agents are uncertain about the true probability
measure. Their investors maximize their utility with respect to
nonadditive probability measures. It turns out that uncertainty
aversion leads to inertia: the agents do not trade the asset
unless the price exceeds or falls below a certain threshold.

We provide a mathematical framework for modelling investor inertia
in a simple microstructure model where asset prices result from
the demand of a large number of small investors whose trading
behavior exhibits inertia. To each agent $a$, we associate a
stationary semi-Markov process $x^a=(x^a_t)_{t \geq 0}$ on a
finite state space which represents the agent's propensity for
trading. The processes $x^a$ have heavy-tailed sojourn times in
some designated ``inert'' state, and relatively thin-tailed
sojourn times in various other states. Semi-Markov processes are
tailor made to model individual traders' inertia as they
generalize Markov processes by removing the requirement of
exponentially distributed, and therefore thin-tailed, holding
times. In addition, we allow for a market-wide amplitude process
$\Psi$, that describes the evolution of typical trading {\em size}
in the market. It is large on heavy-trading days and small on
light trading days. We adopt a non-Walrasian approach to asset
pricing and assume that prices move in the direction of market
imbalance. We show that in a model with many inert investors, long
range dependence in the price process emerges.



\subsection{Long Range Dependence in Financial Time Series}

The observation of long range dependence (sometimes called the
Joseph effect) in financial time series motivated the use of
fractional Brownian motion as a basis for asset pricing models;
see, for instance, \cite{mandelbrot} or \cite{kopp}. By our
invariance principle, the drift-adjusted logarithmic price process
converges weakly to a stochastic integral with respect to a
fractional Brownian motion with Hurst coefficient $H >
\frac{1}{2}$. Our approach may thus be viewed as a microeconomic
foundation for these models. A recent paper that proposes entirely
different economic foundations for models based on \fbm is
\cite{Klupp}. An approximation result for \fbm in the context of a
binary market model is given in \cite{Sottinen}.

As is well known, \fbm processes are not semimartingales, and so
these models may theoretically allow arbitrage opportunities.
Explicit arbitrage strategies for various models were constructed
in \cite{Rogers}, \cite{cheridito} and \cite{erhan}. These
strategies capitalize on the smoothness of \fbm (relative to
standard Brownian motion) and involve rapid trading to exploit the
fine-scale properties of the process' trajectories. As a result,
in our microstructure model, arbitrage opportunities may arise for
other, sufficiently sophisticated, market participants who are
able to take advantage of inert investors by trading frequently.
We discuss a simple combination of both inert and active traders
in Section \ref{twothree}.

\begin{figure}[htb]
  \centering
  \includegraphics[totalheight=4in]{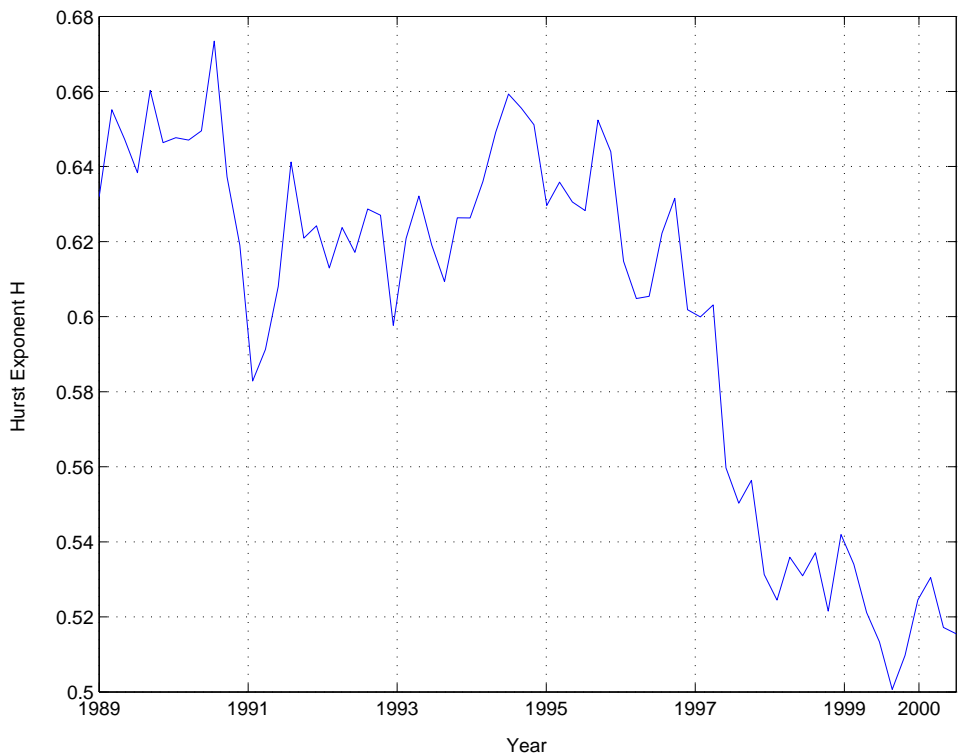}
  \caption{\small{\em Estimates of the Hurst exponent of the S\&P 500 index
  over 1990s, taken from Bayraktar, Poor and Sircar (\cite{bps}).}}
  \label{Hests}
\end{figure}

Evidence of long-range dependence in financial data is discussed
in \cite{kopp}. Bayraktar {\em et al.} \cite{bps} studied an
asymptotically efficient wavelet-based estimator for the Hurst
parameter, and analyzed high frequency S\&P 500 index data over
the span of 11.5 years (1989-2000). It was observed that, although
the Hurst parameter was significantly above the efficient markets
value of $H=\frac{1}{2}$ up through the mid-1990s, it started to
fall to that level over the period 1997-2000 (see Figure
\ref{Hests}). They suggested that this behavior of the market
might be related to the increase in Internet trading, which is
documented, for example, in NYSE's Stockownership2000 \cite{NYSE},
\cite{odean1}, and \cite{choi}, who find that ``after 18 months of
access, the Web effect is very large: trading frequency doubles.''
Barber and Odean \cite{odean2} find that ``after going online,
investors trade more actively, more speculatively and less
profitably than before". Similar empirical findings were recently
reached, using a completely different statistical technique in
\cite{bianchi}. Thus, the dramatic fall in the estimated Hurst
parameter in the late 1990s can be thought of as {\em a
posteriori} validation of the link our model provides between
investor inertia and long-range dependence in stock prices.

We note the evidence of long memory in stock price returns is
mixed. There are several papers in the empirical finance
literature providing evidence for the existence of long memory,
yet there are several other papers that contradict these empirical
findings; see e.g. \cite{bps} for an exposition of this debate and
references. However, long memory is a well accepted feature in
volatility (squared and absolute returns) and trading volume (see
e.g. \cite{contstylized} and \cite{engle}). The mathematical
results of this paper might also be seen as an intermediate step
towards a microstructural foundation for this phenomenon.

\subsection{Mathematical Contributions}

We establish a functional central limit theorem for semi-Markov
processes (Theorem \ref{thm1} below) which extends the results of
Taqqu et. al. \cite{taqqu}, who proved a result similar to ours
for on/off processes, that is, semi-Markov processes taking values
in the binary state space $\{0,1\}$. Their arguments do not carry
over to models with more general state spaces. Our
approach builds on Markov renewal theory. 
We also demonstrate (see Example \ref{exampleA}) that there may be
a different limit behavior when the semi-Markov processes are
centered, a situation which cannot arise in the binary case. Taqqu
and Levy \cite{taqqu2} considered renewal reward processes with
heavy tailed renewal periods and independent and identically
distributed rewards. They assume a general state space, but the
distributions of the length of renewal periods does not depend on
the current state; for an extension to the case of heavy-tailed
rewards, see \cite{taqqu3}. A recent paper \cite{Mikosch} studies
the binary case under a different limit taking mechanism; see also
\cite{Gaigalas-Kaj}.

Binary state spaces are natural for modelling internet traffic,
but for many applications in Economics or Queueing Theory, it is
clearly desirable to have more flexible results that apply to
general semi-Markov processes on finite state spaces. In the
context of a financial market model, it is natural to allow for
both positive (buying), negative (selling) and a zero (inactive)
state. Our results also have applications to complex multi-level
queueing networks where the level-dependent holding-time
distributions are allowed to have slowly decaying tails. They may
serve as a mathematical basis for proving heavy-traffic limits in
the network models studied in, e.g. \cite{Duffield-Whitt},
\cite{D-W} and \cite{Whitt}.

We allow for limits which are integrals with respect to fractional
Brownian motion proving an approximation result for stochastic
integrals of continuous semimartingales with respect to fractional
Brownian motion. Specifically, we consider a sequence of good
semimartingales $\{\Psi^n\}$ and a sequence of stochastic
processes $\{X^n\}$ having zero quadratic variation and give
sufficient conditions which guarantee that joint convergence of
$(X^n,\Psi^n)$ to $(B^H,\Psi)$, where $B^H$ is a \fbm process with
Hurst parameter $H > \frac{1}{2}$, and $\Psi$ is a continuous
semimartingale, implies the convergence of the stochastic
integrals $\int \Psi^n dX^n$ to $\int \Psi d B^H$. In addition, we
obtain a stability result for the integral of a \fbm with respect
to itself. These results may be viewed as an extension of Theorem
2.2 in \cite{Kurtz-Protter} beyond the semimartingale setting.

The remainder of this paper is organized as follows. In Section 2,
we describe the financial market model with inert investors and
state the main result. Section 3 proves a central limit theorem
for stationary semi-Markov processes. Section 4 proves an
approximation result for stochastic integrals of continuous
semimartingales with respect to fractional Brownian motion.


\section{The microeconomic setup and the main results} \label{Section-Model}

We consider a financial market with a set $\ma:=\{a_1,a_2, \ldots
,a_N\}$ of \textsl{agents} trading a single risky asset. Our aim
is to analyze the effects investor inertia has on the dynamics of
stock price processes. For this we choose the simplest possible
setup. In particular, we model right away the behavior of
individual traders rather than characterizing agents' investment
decisions as solutions to individual utility maximization
problems. Such an approach has also been taken in \cite{Garman},
\cite{Foellmer-Schweizer}, \cite{Lux}, \cite{FHK} and \cite{Kruk}
for example.

We associate to each agent $a \in \ma$ a continuous-time
stochastic process $x^a = (x^a_t)_{t \geq 0}$ on a finite state
space $E$, containing zero. This process describes the agent's
\textsl{trading mood}. He accumulates the asset at a rate $\Psi_t
x^a_t$ at time $t \geq 0$. The random quantity $\Psi_t > 0$
describes the size of a typical trade at time $t$, and $x^a_t$ may
be negative, indicating the agent is selling. Agents do not trade
at times when $x^a_t = 0$. We therefore call the state $0$ the
agents' \textsl{inactive state}.

\begin{remark}
In the simplest setting, $x^a\in\{-1,0,1\}$, so that each investor
is either buying, selling or inactive, and $\Psi\equiv 1$: there
is no external amplification. Even here, the existing results in
\cite{taqqu} do not apply because the state space is not binary.
\end{remark}

The holdings of the agent $a \in \ma$ and the ``market imbalance''
at time $t \geq 0$ are given by
\begin{equation}
\label{process-hat-X}
    \int_0^t \Psi_s x^a_s ds \qquad \mbox{and} \qquad
    I_t^N := \sum_{a \in \ma} \int_0^t \Psi_s x^a_s ds,
\end{equation}
respectively. Hence the process $(I^N_t)_{t \geq 0}$ describes the
stochastic evolution of the \textsl{market imbalance}.  In our
microstructure model, market imbalance will be the only component
driving the dynamics of asset prices. 
 All the orders are received
by a single market maker who clears the trades and sets prices as
to reflect the incoming order flows. That is, the market maker
sets prices in reaction to the evolution of market imbalances.

\begin{remark}
In our continuous time model buyers and sellers arrive at
different points in time. Hence the economic paradigm that a
Walrasian auctioneer can set prices such that the markets clear at
the end of each trading period does not apply. Rather, temporary
imbalances between demand and supply will occur, and prices are
assumed to reflect the extent of the current market imbalance. In
the terminology explained in \cite{Garman}, ours is a model of a
``continuous market (trading asynchronously during continuous
intervals of time)", rather than a ``call market (trading
synchronously at pre-established discrete times)". As Garman
reports, the New York Stock Exchange was a call market until 1871,
and since then has become a continuous market. (See also Chapter 1
of \cite{ohara}.)
\end{remark}

We consider the pricing rule
\begin{equation}
\label{process-S}
    dS^N_t = \sum_{a \in \ma} \Psi_t x^a_t dt \quad \mbox{and so} \quad
    S^N_t = S_0 + I_t^N,
\end{equation}
for the evolution of the logarithmic stock price process
$S^N=(S_t^N)_{t \geq 0}$. This is the simplest mechanism by which
incoming buy orders increase the price and sell orders decrease
the price. Other choices might be utilized in a future work
studying, for example, the effect of a nonlinear market depth
function, but these are beyond the scope of the present work. (The
choice of modelling the log-stock price is simply standard finance
practice to define a positive price).

Kruk (\cite{Kruk}) considered a model for continuous auction
market, in which order arrivals are modelled by independent
renewal processes. There are a finite number of possible prices,
and agents randomly submit price dependent limit orders. These are
stored in the order book waiting the arrival of matching orders.
Kruk finds a limiting distribution of the outstanding number of
buy/sell orders at one of the possible prices. In contrast, our
aim is to find the limiting price process that is driven by the
market imbalance under different assumptions on the market
micro-structure.

%



\subsection{The dynamics of individual behavior}

Next, we specify the probabilistic structure of the processes
$x^a$. We assume that the agents are homogeneous and that all the
processes $x^a$ and $\Psi$ are independent. It is therefore enough
to specify the dynamics of some reference process $x = (x_t)_{t
\geq 0}$. In order to incorporate the idea of market inertia as
defined by Assumption \ref{Assumption-tails} below, we assume that
$x$ is a \textsl{semi-Markov} process defined on some probability
space $(\Omega,\mf,\Prob)$ with a finite state space $E$. Here $E$
may contain both positive and negative values and we assume $0 \in
E$. The process $x$ is specified in terms of random variables
$\xi_n: \Omega \rightarrow E$ and $T_n: \Omega \rightarrow \re_+$
which satisfy $0 = T_0 \leq T_1 \leq \cdots$ almost surely and
\[
    \Prob \{\xi_{n+1} = j, T_{n+1}-T_{n} \leq t \big|
    \xi_{1},...,\xi_{n}; T_{1},...,T_{n} \} =
    \Prob \{\xi_{n+1} = j, T_{n+1}-T_{n}\leq t \big|\xi_{n}\}
\]
for each $n \in \n$, $j \in E$ and all $t \in \re_+$ through the
relation
\begin{equation}
    x_t = \sum_{n \geq 0} \xi_n 
    \ind_{[T_n,T_{n+1})}(t).
\end{equation}

\begin{remark}
In economic terms, the representative agent's mood in the random
time interval $[T_{n}, T_{n+1})$ is given by $\xi_n$. The
distribution of the length of the interval $T_{n+1} - T_n$ may
depend on the sequence $\{\xi_n\}_{n \in \n}$ through the states
$\xi_{n}$ and $\xi_{n+1}$. This allows us to assume different
distributions for the lengths of the agents' active and inactive
periods, and in particular to model inertia as a heavy-tailed
sojourn time in the zero state.
\end{remark}

\begin{remark}\label{twofour}
In the present analysis of investor inertia, we do not allow for
feedback effects of prices into agents' investment decisions.
While such an assumption might be justified for small,
non-professional investors, it is clearly desirable to allow
active traders' investment decisions to be influenced by asset
prices. 
When such feedback effects are allowed, the analysis of the price
process is typically confined to numerical simulations because
such models are difficult to analyze on an analytical level.
%
%
An exception is a recent paper \cite{FHK} where the impact of
contagion effects on the asymptotics of stock prices is analyzed
in a mathematically rigorous manner. One could also consider the
present model as applying to (Internet or new economy) stocks
where no accurate information about the actual underlying
fundamental value is available. In such a situation, price is not
always a good indicator of value and is often ignored by
uninformed small investors.
%
\end{remark}

We assume that $x$ is temporally homogeneous under the measure
$\Prob$, that is,
\begin{equation}
\label{time-homogeneous}
    \Prob \{\xi_{n+1} = j, ~T_{n+1}-T_{n} \leq t \big| \xi_{n} = i\} =
    Q(i,j,t)
\end{equation}
is independent of $n \in \n$. By Proposition 1.6 in
\cite{Cinlar75}, this implies that $\{\xi_n\}_{n \in \n}$ is a
homogeneous Markov chain on $E$ whose transition probability
matrix $P=(p_{ij})$ is given by
\[
    p_{ij} = \lim_{t \rightarrow \infty}Q(i,j,t).
\]
Clearly, $x$ is an ordinary temporally homogeneous Markov process
if $Q$ takes the form
\begin{equation} \label{Markov-kernel}
    Q(i,j,t) = p_{ij} \left(1 - e^{-\lambda_i t} \right).
\end{equation}
We assume that the \textsl{embedded Markov chain} $\{\xi_n\}_{n
\in \n}$ satisfies the following condition.

\begin{assumption}
\label{Assumption-pi} For all $i,j \in E$, $i \neq j$ we have that
$p_{ij} > 0$. In particular, there exists a unique probability
measure $\pi$ on $E$ such that $\pi P = \pi$.
\end{assumption}

The conditional distribution function of the length of the $n$-th
sojourn time, $T_{n+1} - T_n$, given $\xi_{n+1}$ and $\xi_n$ is
specified in terms of the \textsl{semi-Markov kernel} $\{
Q(i,j,t); i,j \in E, ~ t \geq 0\}$ and the transition matrix $P$
by
\begin{equation} \label{def-G}
    G(i,j,t) := \frac{Q(i,j,t)}{p_{ij}} = \Prob\{T_{n+1} - T_n \leq t |
    \xi_{n} = i, ~ \xi_{n+1} = j\}.
\end{equation}

For later reference we also introduce the distribution of the
first occurrence of state $j$ under $\Prob$, given $x_0 = i$.
Specifically, for $i\neq j$, we put
\begin{equation} \label{def-travel}
    F(i,j,t) := \Prob\{ \tau_j \leq t | x_0 = i \},
\end{equation}
where $\tau_{j} := \inf \{t \geq 0 : x_t = j\}$. We denote by
$F(j,j,\cdot)$ the distribution of the time until the next
entrance into state $j$ and by
\begin{equation}\label{etamj}
\eta_j := \int t F(j,j,dt)
\end{equation}
the expected time between two occurrences of state $j \in E$.
Further, we recall that a function $L: \re_+ \rightarrow \re_+$ is
called \textsl{slowly varying at infinity} if
\[
    \lim_{t \rightarrow \infty} \frac{L(xt)}{L(t)} = 1 \qquad \mbox{for
    all} \qquad x>0
\]
and that $f(t) \sim g(t)$ for two functions $f,g: \re_+
\rightarrow \re_+$ means $\lim_{t\rightarrow
\infty}\frac{f(t)}{g(t)}=1$.

\begin{assumption} \label{Assumption-tails}
\begin{rmenumerate}
    \item The average sojourn time at state $i \in E$ is finite:
\begin{equation}
\label{mean-time}
    m_i := \E[T_{n+1} - T_n | \xi_n=i] < \infty.
\end{equation}
    Here $\E$ denotes the expectation operator with respect to
    $\Prob$.
    \item There exists a constant $1 < \alpha < 2$ and a locally bounded function $L: \re_+
    \rightarrow \re_+$ which is slowly varying at infinity such that
\begin{equation}
\label{fat-tail}
    \Prob \{ T_{n+1} - T_{n} \geq t \big| \xi_{n} = 0 \}
    \sim t^{-\alpha}L(t).
\end{equation}
    \item The distributions of the sojourn times at state $i \neq 0$ satisfy
\[
    \lim_{t \rightarrow \infty} \frac{\Prob\{T_{n+1} - T_{n} \geq t \big| \xi_{n} = i\}}
    {t^{-(\alpha+1)} L(t)} = 0.
\]
    \item The distribution of the sojourn times in the various
    states have continuous and bounded densities with respect to Lebesgue measure
    on $\re_+$.
\end{rmenumerate}
\end{assumption}

Our condition (\ref{fat-tail}) is satisfied if, for instance, the
length of the sojourn time at state $0 \in E$ is distributed
according to a Pareto distribution. Assumption
\ref{Assumption-tails} (iii) reflects the idea of market inertia:
the probability of long uninterrupted trading periods is small
compared to the probability of an individual agent being inactive
for a long time. It is stronger than the corresponding assumption
for the binary case in \cite{taqqu} where the sojourn time in the
only other state may in fact be as heavy tailed. For our economic
application, however, it is natural to think of the sojourn times
in the various active states as being  thin tailed, such as in the
exponential distribution, since small investors typically do not
trade continually for long periods. 


\subsection{An invariance principle for semi-Markov processes}

In this section, we state our main results. With our choice of
scaling, the logarithmic price process can be approximated in law
by the stochastic integral of $\Psi$ with respect to fractional
Brownian motion $B^H$ where the \textsl{Hurst coefficient} $H$
depends on $\alpha$. The convergence concept we use is weak
convergence on the Skorohod space $\mathbb{D}$ of all real-valued
right continuous processes with left limits. We write
$\mathcal{L}\mbox{-}\lim_{n \rightarrow \infty} Z^n = Z$ if a
sequence of $\mathbb{D}$-valued stochastic processes $\{Z^n\}_{n
\in \mathbb{N}}$, converges in distribution to $Z$.

In order to derive our approximation result, we assume that the
semi-Markov process $x$ is stationary. Under Assumption
\ref{Assumption-pi}, stationarity can be achieved by a suitable
specification of the common distribution of the initial state
$\xi_0$ and the initial sojourn time $T_1$. We denote the
distribution of the stationary semi-Markov processes by $\Prob^*$.
The proof follows from Theorem 4.2.5 in \cite{BFL}, for example. 

\begin{lemma} \label{lemma-stationary}
In the stationary setting, that is, under the law $\Prob^*$ the
following holds:
\begin{rmenumerate}
\item The joint distribution of the initial state and the initial
sojourn time takes the form
\begin{equation}
    \Prob^* \left\{ \xi_0 = k, T_1 > t \right\}
    = \frac{\pi_{k}}{\sum_{j \in E} \pi_{j}m_j}\int_{t}^{\infty}h(k,s)
    ds.\label{cross}
\end{equation}
Here $m_i$ denotes the mean sojourn time in state $i \in E$ as
defined by (\ref{mean-time}), and for $i \in E$,
\begin{equation}
    h(i,t)=1 - \sum_{j \in E} Q(i,j,t)\label{hdef}
\end{equation}
is the probability that the sojourn time at state $i \in E$ is
greater than $t$. \item The law $\nu = (\nu_k)_{k \in E}$ of $x_t$
in the stationary regime is given by
\begin{equation} \label{nu}
    \nu_k = \frac{\pi_k m_k}{\sum_{j \in E} \pi_j m_j}.
\end{equation}
\item The conditional joint distribution of $(\xi_1, T_1)$, given
$\xi_0$ is
\begin{equation}
    \Prob^* \left\{ \xi_1 = j, T_1 < t\mid \xi_0 = k \right\}
    = \frac{p_{kj}}{m_{k,j} } \int_0^t [ 1 - G(k,j,s) ]
    ds.\label{star}
\end{equation}
Here $m_{k,j} := \int_0^\infty [1-G(k,j,s)]ds$ denotes the
conditional expected sojourn time at state $k$, given the next
state is $j$, and the functions $G(k,j,\cdot)$ are defined in
(\ref{def-G}).
\end{rmenumerate}
\end{lemma}


Let us now introduce a dimensionless parameter $\eps>0$, and
consider the rescaled processes $x^a_{t/\eps}$. For $\eps$ small,
$x^a_{t/\eps}$ is a ``speeded-up" semi-Markov process. In other
words, the investors' individual trading dispensations are
evolving on a faster scale than $\Psi$. Observe, however, that we
are not altering the main qualitative feature of the model. That
is, agents still remain in the inactive state for relatively much
longer times than in an active state.

Mathematically, there is no reason to restrict ourselves to the
case where $\Psi$ is non-negative. Hence we shall from now on only
assume that $\Psi$ is a continuous semimartingale. Given the
processes $\Psi$ and $x^a$ $(a \in \{a_1, \ldots ,a_N\})$, the
aggregate order rate at time $t$ is given by
\begin{equation}
    Y^{\eps,N}_t = \sum_{a \in \ma} \Psi_t x^a_{t/\eps} \label{YepsN}.
\end{equation}
Let $\mu := \E^* x_t$ and $X^{\eps,N} = (X^{\eps,N}_t)_{0 \leq t
\leq T}$ $(T > 0)$ be the centered process defined by
\begin{equation}
    X^{\eps,N}_t
    :=\int_0^t\sum_{a\in\ma}\Psi_s(x^a_{s/\eps}-\mu)\,ds=
    \int_0^{t} Y^{\eps,N}_s\,ds - \mu N\int_{0}^{t}
    \Psi_{s}ds.
\label{defXepsN}
\end{equation}
We are now ready to state our main result. Its proof will be
carried out in Sections \ref{Section-invariance} and
\ref{Section-sigma}. The definition of the stochastic integral
with respect to \fbm will also be given in Section
\ref{Section-sigma}.

%

\begin{theorem}
\label{thm1} Let $\Psi=(\Psi_t)_{t \geq 0}$ be a continuous
semimartingale on $(\Omega,\mf,\Prob^*)$ with a decomposition
$\Psi= M+A$, in which $M$ is a local martingale and $A$ is an
adapted process of finite variation. We assume that
$\mathbb{E}\{[M,M]_{T}\}<\infty$ and that
$\mathbb{E}\{|A|_{T}\}<\infty$, where $(|A|_{t})_{t \ge 0}$ is the
total variation of $A$. If Assumptions \ref{Assumption-pi} and
\ref{Assumption-tails} are satisfied, and if $\mu \sum_{k \in E} k
\frac{m_k}{\eta_k^2}
> 0$, then there exists $c > 0$ such that the process $X^{\eps,N}$
satisfies
\begin{equation}
\label{limit}
    \mathcal{L}\mbox{-}\lim_{\eps\downarrow 0} \mathcal{L}\mbox{-}\lim_{N \rightarrow
    \infty} \left(\frac{1}{\eps^{1-H}\sqrt{NL(\eps^{-1})}}X^{\eps,N}_t
    \right)_{0 \leq t \leq T} = \left(c \int_0^t \Psi_s d B^H_s \right)_{0 \leq t \leq T}.
\end{equation}
Here the Hurst coefficient of the fractional Brownian motion
process $B^H$ is $H = \frac{3-\alpha}{2} > \frac{1}{2}$.
\end{theorem}

Observe that Theorem \ref{thm1} does not apply to the case $\mu =
0$. For centered semi-Markov processes $x^a$, Example
\ref{exampleA} below illustrates that the limiting process depends
on the tail structure of the waiting time distribution in the
various active states. This phenomenon does not arise in the case
of binary state spaces.

\begin{remark} \begin{rmenumerate}
\item Theorem \ref{thm1} says the drift-adjusted logarithmic price
process in our model of inert investors can be approximated in law
by the stochastic integral of $\Psi$ with respect to a fractional
Brownian motion process with Hurst coefficient $H > \frac{1}{2}$.


\item In a situation where the processes $x^a$ are independent,
stationary and ergodic Markov processes on $E$, that is, in cases
where the semi-Markov kernel takes the form (\ref{Markov-kernel}),
it is easy to show that
\[
    \mathcal{L}\mbox{-}\lim_{\eps\downarrow 0}
    \mathcal{L}\mbox{-}\lim_{N \rightarrow \infty}
    \left(\frac{1}{\sqrt{\eps N}}X^{\eps,N}_t
    \right)_{0 \leq t \leq T} = \left( c \int_0^t \Psi_s d W_s
    \right)_{0 \leq t \leq T}
\]
where $(W_t)_{t \geq 0}$ is a standard Wiener process. Thus, if
the market participants are not inert, that is, if the
distribution of the lengths of the agents' inactivity periods is
thin-tailed, no arbitrage opportunities emerge because the limit
process is a
semimartingale. 
\end{rmenumerate}
\end{remark}

The proof of Theorem \ref{thm1} will be carried out in two steps.
In Section \ref{Section-invariance} we prove a functional central
limit theorem for stationary semi-Markov processes on finite state
spaces. In Section \ref{Section-sigma} we combine our central
limit theorem for semi-Markov processes with extensions of
arguments given in \cite{Kurtz-Protter} to obtain (\ref{limit}).

\subsection{Markets with both Active and Inert Investors\label{twothree}}
It is simple to extend the previous analysis to incorporate both
active and inert investors. Let $\rho$ be the ratio of active to
inert investors. We associate to each \textsl{active} trader $b
\in \{1,2, \ldots, \rho N\}$ a stationary Markov chain $y^b =
(y^b_t)_{t \geq 0}$ on the state space $E$. The processes $y^b$
are independent and identically distributed and independent of the
processes $x^a$. The thin-tailed sojourn time in the zero state of
$y^b$ reflects the idea that, as opposed to inert investors, these
agents frequently trade the stock. We assume for simplicity that
$\Psi \equiv 1$. With $\hat{Y}^{\eps,N}_t = \sum_{b = 1}^{\rho N}
\left(y^b_{t/\eps} - \E^* y_0 \right)$ and $\hat{X}^{\eps,N}_t :=
\int_0^t \hat{Y}^{\eps,N}_s ds$, it is straightforward to prove
the following modification of Theorem \ref{thm1}.

\begin{proposition}
    Let $x^a$ $(a = 1,2, \ldots ,N)$ be semi-Markov processes that satisfy the assumption
    of Theorem \ref{thm1}. If $y^b$ $(b=1,2, \ldots ,\rho N)$ are independent
    stationary Markov processes on $E$, then there exist
    constants $c_1,c_2 > 0$ such that
\[
    \mathcal{L}\mbox{-}\lim_{\eps\downarrow 0} \mathcal{L}\mbox{-}\lim_{N \rightarrow
    \infty}
    \left(\frac{1}{\eps^{1-H}\sqrt{NL(\eps^{-1})}}X^{\eps,N}_t +
    \frac{1}{\sqrt{N \eps}}\hat{X}^{\eps,N}_t
    \right)_{0 \leq t \leq T} = \left(c_1 B^H_t + c_2\sqrt{\rho}\, W_t \right)_{0 \leq t \leq
    T}.
\]
    Here, $W=(W_t)_{t \geq 0}$ is a standard Wiener process.
\end{proposition}

Thus, in a financial market with both active and inert investors,
the dynamics of the asset price process can be approximated in law
by a stochastic integral with respect to a superposition, $B^H +
\delta W$, of a fractional and a regular Brownian motion. It is
known (\cite{cheridito-mixed}) that $B^H + \delta W$ is a
semimartingale for any $\delta\neq 0$, if $H> \frac{3}{4}$,  that
is, if $\alpha < \frac{3}{2}$, but not if $H \in
(\frac{1}{2},\frac{3}{4}]$. Thus, no arbitrage opportunities arise
if the small investors are ``sufficiently inert.'' The parameter
$\alpha$ can also be viewed as a measure for the fraction of small
investors that are active at any point in time. Hence, independent
of the actual trading volume, the market is arbitrage free in
periods where the fraction of inert investors who are active on
the financial market is small enough.


\section{A limit theorem for semi-Markov processes}
\label{Section-invariance}

This section establishes Theorem \ref{thm1} for the special case
$\Psi\equiv 1$. We approach the general case where $\Psi$ is a
continuous semimartingale in Section \ref{Section-sigma}. Here we
consider the situation where
\[
    Y^{\eps,N}_t = \sum_{a \in \ma} x^a_{t/\eps}
    \quad  \mbox{and where} \quad X^{\eps,N}_t =
    \int_0^tY^{\eps,N}_s\,ds - N\mu t,
\]
and prove a functional central limit theorem for stationary
semi-Markov processes. Our Theorem \ref{thm2} below extends the
results in \cite{taqqu} to situations where the semi-Markov
process takes values in an arbitrary finite state space. The
arguments given there are based on results from ordinary renewal
theory, and do not carry over to models with more general state
spaces. The proof of the following theorem will be carried out
through a series of lemmas.

\begin{theorem} \label{thm2}
Let $H = \frac{3-\alpha}{2}$. Under the assumptions of Theorem
\ref{thm1},
\begin{equation} \label{limit1}
    \mathcal{L}\mbox{-}\lim_{\eps\downarrow 0} \mathcal{L}\mbox{-}\lim_{N
    \rightarrow \infty} \left(\frac{1}{\eps^{1-H}\sqrt{NL(\eps^{-1})}}X^{\eps,N}_t
    \right)_{0 \leq t \leq T} =  \left(c B^H_t \right)_{0 \leq t \leq T}.
\end{equation}
\end{theorem}

Let $\gamma$ be the covariance function of the semi-Markov process
$(x_t)_{t \geq 0}$ under $\Prob^*$, and consider the case
$\eps=1$. By the Central Limit Theorem, and because $x$ is
stationary, the process $Y=(Y_t)_{t\geq 0}$ defined by
\begin{equation} \label{def-Y}
    Y_t = \mathcal{L}\mbox{-}\lim_{N \rightarrow \infty}
    \frac{1}{\sqrt{N}}(Y^{1,N}_t-N\mu)
\end{equation}
is a stationary zero-mean Gaussian process. It is easily checked
that the covariance function of the process
$(\frac{1}{\sqrt{N}}Y^{1,N}_t)$ is also $\gamma$ for any $N$, and
hence for $Y_t$. By standard calculations, the variance of the
aggregate process $(\int_0^tY_s\,ds)$ at time $t \geq 0$ is given
by
\begin{equation}
    \Var(t) := \Var\left(\int_0^tY_s\,ds\right) = 2 \int_0^t
    \left( \int_0^v \gamma(u) du \right) dv.
    \label{representation-variance}
\end{equation}

In the first step towards the proof of Theorem \ref{thm2}, we can
proceed by analogy with \cite{taqqu}. We are interested in the
asymptotics as $\eps\downarrow 0$ of the process
\begin{equation}
    X^\eps_t := \int_0^t Y_{s/\eps}\,ds, \label{Xepsdef}
\end{equation}
which can be written $X^\eps_t = \eps\int_0^{t/\eps}Y_s\,ds$.
Therefore the object of interest is the large $t$ behavior of
$\Var(t)$. Suppose that we can show
\begin{equation}
\label{variance}
    \Var(t) \sim c^2 t^{2H} L(t) \quad \mbox{as} \quad
    t \rightarrow \infty.
\end{equation}
Then the mean-zero Gaussian processes $X^\eps=(X^\eps_t)_{t \geq
0}$ have stationary increments and satisfy
\begin{equation} \label{second-moment}
    \lim_{\eps\downarrow 0} \E^*\left(\frac{1}{\eps^{1-H}
    \sqrt{L(\eps^{-1})}}X^\eps_t\right)^2 = c^2 t^{2H}.
\end{equation}
Since the variance characterizes the finite dimensional
distributions of a mean-zero Gaussian process with stationary
increments, we see that the finite dimensional distributions of
the process $\left( \frac{1}{ \eps^{1-H}\sqrt{L(\eps^{-1})}}
X^{\eps}_{t} \right)_{t \geq 0}$ converge to $\left(c
B^{H}_{t}\right)_{t \geq 0}$ whenever (\ref{variance}) holds. The
following lemma gives a sufficient condition for (\ref{variance})
in terms of the covariance function $\gamma$.

\begin{lemma} \label{lemma-covariance}
For (\ref{variance}) to hold, it suffices that
\begin{equation}
\label{variance-gamma}
    \gamma(t) \sim c^2 H(2H-1)t^{2H-2}L(t) \quad \mbox{as} \quad t
    \rightarrow \infty.
\end{equation}
\end{lemma}
\begin{proof}
By Proposition 1.5.8 in \cite{bingham}, every slowly varying
function $L$ which is locally bounded on $\re_+$ satisfies
\[
    \int_0^t \tau^\beta L(\tau) d \tau \sim  \frac{t^{\beta+1} L(t)}{\beta +
    1}
\]
if $\beta > -1$. Applying this proposition to the slowly varying
function
\[
     \tilde{L}(t):=\frac{\gamma(t)}{c^{2}H(2H-1)t^{2H-2}},
\]
we conclude
\[
    \int_{0}^{t}\int_{0}^{v}\gamma(u) du\,dv \sim \frac{c^{2}}{2} t^{2H} L(t),
\]
and so our assertion follows from (\ref{representation-variance}).
\end{proof}

Before we proceed with the proof of our main result, let us
briefly consider the case $\mu = 0$ which is not covered by our
theorem. For semi-Markov processes whose ``heavy-tailed state''
happens to be the mean, the structure of the limit process depends
on the distribution of the sojourn times in the various active
states.\footnote{We thank Chris Rogers for Example
\ref{exampleA}.}

\begin{example} \label{exampleA} We consider the case
$E=\{-1,0,1\}$, and assume that $p_{-1,0}=p_{1,0} = 1$ and that
$p_{0,-1}=p_{0,1} = \frac{1}{2}$. With $\nu_1 = \Prob^*\{x_t=1\} >
0$, we obtain
\[
    \gamma(t) = \nu_1 \left( \E^*[x_t x_0 | x_0 = 1] + \E^*[x_t x_0
    | x_0 = -1] \right).
\]

Suppose that the sojourns in the inactive state are heavy tailed,
and that the waiting times in the active states are exponentially
distributed with parameter 1. In such a symmetric situation
\[
    \E^*[x_t x_0 | x_0 = \pm  1] = \Prob^* \{T_1 \geq t | x_0 = \pm 1\} =
    e^{-t}.
\]
Therefore, $\gamma(t) =  2 \nu_1 e^{-t}$. In view of
(\ref{second-moment}), this yields $c > 0$ such that
\[
    \mathcal{L}\mbox{-}\lim_{\eps\downarrow 0}
    \mathcal{L}\mbox{-}\lim_{N \rightarrow \infty}
    \left(\frac{1}{\sqrt{\eps N}}X^{\eps,N}_t \right)_{0 \leq t \leq T}
    = \left( c W_t \right)_{0 \leq t \leq T}
\]
    for some standard Wiener process $W$.

\end{example}

In order to prove Theorem \ref{thm2}, we need to establish
(\ref{variance-gamma}). For this, the following representation of
the covariance function turns out to be useful: in terms of the
marginal distribution $\nu_i = \Prob^* \{x_t=i\}$ $(i \in E)$ of
the stationary semi-Markov process given in Lemma
\ref{lemma-stationary} (i), and in terms of the conditional
probabilities
\[
    P^{*}_{t}(i,j) := \Prob^*\{x_t = j | x_0 = i\},
\]
we have
\begin{equation}\label{eq:cov}
    \gamma(t) = \sum_{i,j \in E} ij \nu_i \left( P^{*}_{t}(i,j) -
    \nu_j \right).
\end{equation}
It follows from Proposition 6.12 in \cite{Cinlar75}, for example,
that $P^{*}_{t}(i,j) \rightarrow \nu_j$ as $t \rightarrow \infty$.
Hence $\lim_{t \rightarrow \infty} \gamma(t) = 0$. In order to
prove Theorem \ref{thm2}, however, we also need to show that this
convergence is sufficiently slow. We shall see that the agents'
inertia accounts for the slow decay of correlations. It is thus
the agents' inactivity that is responsible for that fact that the
logarithmic price process is not approximated by a stochastic
integral with respect to a Wiener process, but by an integral with
respect to fractional Brownian motion.
%
%

We are now going to determine the rate of convergence of the
covariance function to 0. To this end, we show that
$P^{*}_{t}(i,j)$ can be written as a convolution of a renewal
function with a slowly decaying function plus a term which has
asymptotically, i.e.,  for $t \rightarrow \infty$, a vanishing
effect compared to the first term; see Lemma
\ref{representation-P*} below. We will then apply results from
\cite{heath} and \cite{jelen} to analyze the tail structure of the
convolution term.

Let
\[
    R(i,j,t) := \E\left\{\sum_{n=0}^{\infty}\ind_{\{\xi_{n}=j,
    T_{n} \leq t\}} ~ \big{|} ~ x_0=\xi_0 = i \right\}
\]
be the expected number of visits of the process $(x_t)_{t \geq 0}$
to state $j$ up to time $t$ in the \textsl{non-stationary}
situation, i.e., under the measure $\Prob$, given $x_0 = i$. For
fixed $i,j \in E$, the function $t \mapsto R(i,j,t)$ is a renewal
function. If, under $\Prob$, the initial state is $j$, then the
entrances to $j$ form an ordinary renewal process and
\begin{equation}\label{ordrenewal}
    R(j,j,t) = \sum_{n=0}^{\infty}F^{n}(j,j,t).
\end{equation}
Here $F(j,j,\cdot)$ denotes the distribution of the travel time
between to occurrences of state $j \in E$ as defined in
(\ref{def-travel}), and $F^{n}(j,j,\cdot)$ is the $n$-fold
convolution of $F(j,j,\cdot)$. On the other hand, if $i \neq j$,
the time until the first visit to $j$ has distribution
$F(i,j,\cdot)$ under $\Prob$ which might be different from
$F(j,j,\cdot)$. In this case $R(i,j,\cdot)$ satisfies a delayed
renewal equation, and we have
\begin{equation}
\label{R-non-stationary}
    R(i,j,t)=\int_{0}^{t} R(j,j,t-u) F(i,j,du).
\end{equation}
We refer the interested reader to  \cite{Cinlar75} for a survey on
Markov renewal theory.

Let us now return to the stationary setting and derive a
representation for the expected number $R^*(i,j,t)$ of visits of
the process $(x_t)_{t \geq 0}$ to state $j$ up to time $t$ under
$\Prob^*$, given $x_0 = i$. To this end, we denote by
$F^*(i,j,\cdot)$ the distribution function in the stationary
setting of the first occurrence of $j$, given $x_0 = i$ and put
\[
    P_t(i,j) := \Prob \{x_t = j | x_0=i \}.
\]
Given the first jump time $T_1$ and given that $x_{T_1} = i$ we
have that
\begin{equation}
\label{equal-probabilities}
    \Prob^* \{x_{t} = j | x_{T_1} = i\} = P_{t-T_1}(i,j) \quad
    \mbox{on} \quad \{t \geq T_1\}.
\end{equation}
Thus,
\begin{equation}\label{eq:srij}
    R^{*}(i,j,t)=\int_{0}^{t}R(j,j,t-u) F^{*}(i,j,du).
\end{equation}


\subsection{A representation for the conditional transition
probabilities}

In this section we derive a representation for $P^*_{t}(i,j)$
which will allow us to analyze the asymptotic behavior of
$P^*_{t}(i,j)-\nu_j$. To this end, we recall the definition of the
joint distribution of the initial state and the initial sojourn
time and the definition of the conditional joint distribution of
$(\xi_1, T_1)$, given $\xi_0$ from (\ref{cross}) and (\ref{star})
respectively. We define
\begin{equation}\label{defnofhats}
    s(i,t) := \Prob^*\{\xi_0 = i, T_1 > t\} \quad \mbox{and} \quad
    \hat{s}(i,j,t) = \Prob^* \{\xi_1 = j, T_1 \leq t | \xi_0 = i\}.
\end{equation}
In terms of these quantities, the transition probability
$P^*_t(i,j)$ can be written as
\begin{equation} \label{eq:ueq}
    P^*_{t}(i,j) = \frac{s(i,t)}{\nu_{i}}\delta_{ij}
    + \sum_{k \in E} \int_{0}^{t} P_{t-u}(k,j)
    \hat{s}(i,k,du).
\end{equation}
Here the first term on the right-hand-side of (\ref{eq:ueq})
accounts for the $\Prob^*$-probability that $x_0=i$ and that the
state $i$ survives until time $t$. The quantity $\int_0^t
P_{t-u}(k,j) \hat{s}(i,k,du)$ captures the conditional probability
that the first transition happens to be to state $k$ before time
$t$, given $\xi_0=i$. Observe that we integrate the conditional
probability $P_{t-u}(i,j)$ and not $P^*_{t-u}(i,j)$: conditioned
on the value of semi-Markov process at the first renewal instance
the distributions of $(x_t)_{t \geq 0}$ under $\Prob$ and
$\Prob^*$ are the same; see (\ref{equal-probabilities}).

In the sequel it will be convenient to have the following
convolution operation: let $\tilde{h}$ be a locally bounded
function, and $\tilde{F}$ be a distribution function both of which
are defined on $\re_{+}$. The convolution $\tilde{F}*\tilde{h}$ of
$\tilde{F}$ and $\tilde{h}$ is given by

\begin{equation}\label{conv}
    \tilde{F}*\tilde{h}(t):=\int_{0}^{t}\tilde{h}(t-x)\tilde{F}(dx)
    \quad \mbox{for $t \geq 0$}.
\end{equation}

\begin{remark} \label{rem-renewal}
Since $\tilde{F}*\tilde{h}$ is locally bounded, the map $t \mapsto
G*(\tilde{F}*\tilde{h})(t)$ is well defined for any distribution
$G$ on $\re_+$. Moreover, $G*(\tilde{F}*\tilde{h})(t) =
\tilde{F}*(G*\tilde{h})(t) = (G*\tilde{F})*\tilde{h}(t)$. In this
sense distributions acting on the locally bounded function can
commute. Thus, for the renewal function $R =
\sum_{n=0}^{\infty}\tilde{F}^{n}$ associated to $\tilde{F}$, as
defined in (\ref{ordrenewal}), the integral $R*\tilde{h}(t)$ is
well defined and
$R*(G*\tilde{h})(t)=G*(R*\tilde{h})(t)=(R*G)*\tilde{h}(t)$.
\end{remark}

We are now going to establish an alternative representation for
the conditional probability $P_t^*(i,j)$ that turns out to be more
appropriate for our subsequent analysis.

\begin{lemma} \label{representation-P*}
In terms of the quantities $s(i,t)$ and $h(i,t)$ in (\ref{hdef})
and $R^*(i,j,t)$, we have
\begin{equation}\label{eq:ipij}
    P^*_{t}(i,j)=\frac{s(i,t)}{\nu_{i}}\delta_{ij}+\int_{0}^{t}
    h(j,t-s) R^*(i,j,ds).
\end{equation}
\end{lemma}
\begin{proof} In view of (\ref{eq:ueq}), it is enough to show
\[
    R^*(i,j,t)*h(j,t) = \sum_{k \in E} \int_0^t
    P_{t-u} (k,j) \hat{s}(i,k,du).
\]
To this end, observe first that $F^{*}(i,j,t)$ can be decomposed
as
\begin{equation}\label{F*}
    F^{*}(i,j,t) = \hat{s}(i,j,t) +
    \sum_{k \neq j} \int_{0}^{t} F(k,j,t-u) \hat{s}(i,k,du).
\end{equation}
Indeed, $\hat{s}(i,j,t)$ is the probability that the first
transition takes place before time $t$ and happens to be to state
$j \in E$, and
\[
    \int_{0}^{t} F(k,j,t-u) \hat{s}(i,k,du) = \Prob^* \{x_v = j
    \mbox{ for some $v \leq t$}, ~ x_{T_1} = k | x_0 = i \} .
\]
In view of (\ref{eq:srij}) and (\ref{F*}), Remark
\ref{rem-renewal} yields
\begin{align*}
R^*(i,j,t)*h(j,t) ~
    = & ~ R(j,j,t)*F^*(i,j,t)*h(j,t) ~ \\
    = & ~ R(j,j,t)*\hat{s}(i,j,t)*h(j,t) + \sum_{k \neq j}
    F(k,j,t) * \hat{s}(i,k,t) * R(j,j,t) * h(j,t).
\end{align*}
Now recall from Proposition 6.3 in \cite{Cinlar75}, for example,
that
\[
    P_{t}(i,j)=\int_{0}^{t} h(j,t-s) R(i,j,ds).
\]
Thus, by also using (\ref{R-non-stationary}) we obtain
\begin{eqnarray*}
    R^*(i,j,t)*h(j,t)
    & = & \hat{s}(i,j,t) * P_t(j,j) + \sum_{k \neq j}
    R(j,j,t) * F(k,j,t) * \hat{s}(i,k,t)* h(j,t)  \\
    & = & \hat{s}(i,j,t)*P_t(j,j) + \sum_{k \neq j}
    R(k,j,t) *  \hat{s}(i,k,t) * h(j,t) \\
    & = & \sum _{k} \hat{s}(i,k,t)*P_t(k,j).
\end{eqnarray*}
This proves our assertion.
\end{proof}


\subsection{The rate of convergence to equilibrium}

Now, our goal is to derive the rates of convergence of the
mappings $t \mapsto s(i,t)$ to $0$ and $t \mapsto
R^*(i,j,t)*h(j,t)$ to $\nu_j$, respectively. Due to (\ref{eq:cov})
it is enough to analyze the case $i,j \neq 0$. To this end we
shall first study the asymptotic behavior of the map $t \mapsto
R^*(i,j,t)$. Since $R(j,j,\cdot)$ is a renewal function, we see
from (\ref{eq:srij}) and (\ref{eq:ipij}) that the asymptotics of
$P^*_{t}(i,j)$ can be derived as an application of Theorem
\ref{sidthm} essentially if we can show that
\begin{equation} \label{convergence-F}
   F^*(i,j,t) * h(j,t) = o(\bar{F}(j,j,t)).
\end{equation}

\subsubsection{The tail structure of the travel times}

Let us first deal with the issue of finding the convergence rate
of $\bar{F}(j,j,t)=1-F(j,j,t)$ to 0. To this end, we introduce the
family of random variables
\begin{equation}\label{eq:Y}
     \Theta = \{(\theta_{i,j}^{\ell}), i,j \in E, \ell=0,1,2,...\},
\end{equation}
such that any two random variables in $\Theta$ are independent,
and for fixed pair $(i,j)$ the random variables $\theta_{i,j}^{k}$
have $G(i,j,\cdot)$ as their common distribution function. To ease
the notational complexity we assume that the law of
$T_{n+1}-T_{n}$ only depends on $\xi_{n}$. We shall therefore drop
the second sub-index from the elements of $\Theta$. The random
variables $(\theta_i^\ell)_{\ell \in \mathbb{N}}$ are independent
copies of the sojourn time in state $i$. We shall prove Lemma
\ref{lemma-convergence-Fbar} below under this additional
assumption. The general case where $T_{n+1}-T_{n}$ depends both on
$\xi_{n}$ and $\xi_{n+1}$ can be analyzed by similar means.

Let $N^{i,j}_{k}$ denote the number of times the embedded Markov
chain $\{\xi_{n}\}_{n \in \n}$ visits state $k \in E$ before it
visits state $j$, given $\xi_0= i$. By definition, $N_j^{i,j}=0$
with probability one. We denote by $\vec{N}^{i,j}$ the vector of
length $|E|$ with entries $N^{i,j}_{k}$, and by $\vec{n} =
(n_k)_{k \in E}$ an element of $\mathbb{N}^{|E|}$. Then we have
\begin{equation}\label{eq:F}
     \bar{F}(i,j,t) = \Prob \left\{ \theta_{i}^0 + \sum_{k \neq j}
     \sum_{\ell=1}^{N^{i,j}_{k}} \theta_{k}^{\ell} > t \right\}.
\end{equation}
With $G(k,t) := \Prob\{ T_{n+1}-T_{n} \leq t \big| \xi_{n} = k\}$,
we can rewrite (\ref{eq:F}) as
\begin{equation} \label{eq:F2}
\begin{split}
     \bar{F}(i,j,t) ~ & = ~~ \sum_{\vec{n}}
     \Prob \left\{ \theta_{i}^0 + \sum_{k \neq j} \sum_{\ell=1}^{n_{k}}
     \theta_{k}^{\ell} > t \bigg | \vec{N}^{i,j} = \vec{n}\right\}
     \Prob \{ \vec{N}^{i,j} = \vec{n}\} \\
     ~ & = ~~  1 - G(i,t) * \sum_{\vec{n}}
     \underset{k \neq j}{*} G^{n_{k}}(k,t)
     \Prob \{\vec{N}^{i,j} = \vec{n}\}. 
\end{split}
\end{equation}
Our goal is now to show that
\begin{equation}
\label{convergence-Fjj}
     \lim_{t \rightarrow \infty}
     \frac{\bar{F}(i,j,t)}{t^{-\alpha}L(t)} = \sum_{n \geq 0}
     n \Prob\{N^{i,j}_0 = n \}+\delta_{i,0}
\end{equation}
for $j \neq 0$. Here, $\delta_{i,0}=0$ if $i \neq 0$ and
$\delta_{i,0}=1$ otherwise. The first term on the right-hand-side,
$\sum_{n \geq 0} n \Prob\{N^{i,j}_0 = n \}$, is the expected
number of occurrences of state $0$ under $\Prob$ before the first
visit to state $j$, given $\xi_0 = x_0 = i.$ This quantity is
positive, due to Assumption \ref{Assumption-pi}. In order to prove
(\ref{convergence-Fjj}), we need the following results which
appear as Lemma 10 in \cite{jelen}.

\begin{lemma}\label{domlemma}
Let $F_{1},...,F_{m}$ be probability distribution functions such
that, for all $j \neq i$, we have
$\bar{F}_{j}(t)=o(\bar{F}_{i}(t))$.
Then for any positive integers $n_1, \ldots , n_m$,
\[
     1 - F_{1}^{n_{1}}*...*F_{m}^{n_{m}}(t) \sim n_{i}\bar{F}_{i}(t).
\]
Moreover, for each $u>0$, there exists some $K_{u}<\infty$ such
that
\[
     \frac{1-F_{1}^{n_{1}}*...*F_{m}^{n_{m}}(t)}{1 - F_{i}^{n_i}(t)}
     \leq K_{u}(1+u)^{n_{i}}
\]
for all $t \geq 0$.
\end{lemma}

We are now ready to prove (\ref{convergence-Fjj}).

\begin{lemma} \label{lemma-convergence-Fbar}
Under the assumptions of Theorem \ref{thm2} we have, for $j\neq
0$,
\[
     \lim_{t \rightarrow \infty} \frac{\bar{F}(i,j,t)}{t^{-\alpha}L(t)}
     = \sum_{n \geq 0} n \Prob\{N^{i,j}_0 = n \} + \delta_{i,0} > 0.
\]
\end{lemma}
\begin{proof}
Let us first prove that the expected number of occurrences of
state $0$ before the first return to state $j$ occurs is finite.
To this end, we put $p = \min\{p_{ij}: i,j \in E, i \neq j\} > 0$.
Since
\[
     \Prob\{ N^{i,j}_0=n \} \leq \Prob \{ \xi_m \neq j ~ \mbox{for
     all} ~ m \leq n \} \leq (1-p)^{n},
\]
we obtain
\begin{equation} \label{eq:fv}
     \sum_{n \geq 0} n \Prob\{ N^{i,j}_0=n \} \leq \sum_{n \geq 0}
     n (1-p)^{n} < \infty.
\end{equation}
Now, we define a probability measure $\bar\mu$ on
$\mathbb{N}^{|E|}$ by
\[
     \bar\mu \{\vec{n}\} = \Prob \{\vec{N}^{i,j}=\vec{n}\}
\]
and put
\[
     A_{\vec{n}}(t)= G(i,t) \underset{k \neq
     j}{*}G^{n_{k}}(k,t).
\]
Since $\frac{1-G(0,t)}{t^{-\alpha} L(t)} \rightarrow 1$ as $t
\rightarrow \infty$, the first part of Lemma \ref{domlemma} yields
\[
     \lim_{t \rightarrow \infty}
     \frac{1-A_{\vec{n}}(t)}{t^{-\alpha}L(t)} =
     \lim_{t \rightarrow \infty}\left(\frac{1-G(i,t)\underset{k \neq
     j}{*}G^{n_{k}}(k,t)}{1-G(0,t)}\right)\left(
     \frac{1-G(0,t)}{t^{-\alpha}L(t)}\right) = n_{0} + \delta_{i,0}.
\]
 From the definition of the measure $\bar\mu$, we obtain
\[
     \frac{ 1-\sum_{\vec{n}}
     A_{\vec{n}}(t) \Prob \{\vec{N}^{i,j}=\vec{n}\}}{t^{-\alpha}L(t)}
     = \E_{\bar\mu} \left\{\frac{1-A_{\vec{n}}(t)}{t^{-\alpha}L(t)}\right\},
\]
and so our assertion follows from the dominated convergence
theorem if we can show that
\begin{equation}\label{eq:bn}
     \sup_{t}\frac{1-A_{\vec{n}}(t)}{t^{-\alpha}L(t)} \in L^{1}(\bar{\mu}).
\end{equation}
To verify (\ref{eq:bn}), we will use the second part of Lemma
\ref{domlemma}. For each $u > 0$ there exists a constant $K_{u}$
such that
\[
    \frac{1-A_{\vec{n}}(t)}{t^{-\alpha}L(t)}
    =\left(\frac{1-A_{\vec{n}}(t)}{1-G(0,t)}\right)
    \left(\frac{1-G(0,t)}{t^{-\alpha}L(t)}\right) \leq
    K_{u}(1+u)^{n_{0} + \delta_{i,0}}\sup_{t}\frac{1-G(0,t)}
    {t^{-\alpha}L(t)}.
\]
Since
\[
     \lim_{t\rightarrow \infty}\frac{1-G(0,t)}{t^{-\alpha}L(t)}=1,
\]
and because we are only interested in the asymptotic behavior of
the function $t \mapsto \frac{\bar{F}(i,j,t)}{t^{-\alpha} L(t)}$,
we may with no loss of generality assume that
\[
     \sup_{t}\frac{1-G(0,t)}{t^{-\alpha}L(t)}=1.
\]
This yields
\begin{equation}\label{eq:sup}
     \sup_{t} \frac{1-A_{\vec{n}}(t)}{t^{-\alpha}L(t)}
     \leq K_{u} (1+u)^{n_{0} + \delta_{i,0}}.
\end{equation}
 From (\ref{eq:fv}) and (\ref{eq:sup}) we get
\[
     \E_{\bar\mu} \left\{ \sup_{t}\frac{1-A_{\vec{n}}(t)}{t^{-\alpha}L(t)}
     \right\} \leq K_{u} (1 + u)^{\delta_{i,0}}
     \sum_{k = 0}^{\infty}(1-p)^{k} (1+u)^{k}.
\]
Choosing $u < \frac{p}{1-p}$ we obtain $\beta := (1-p)(1+u) < 1$
and so the assertion follows from
\[
     \E_{\bar\mu}\left\{\sup_{t}\frac{1-A_{\vec{n}}(t)} {t^{-\alpha}L(t)}
     \right\} \leq K_{u} (1 + u)^{\delta_{i,0}}
     \sum_{n \geq 0} \beta^{n} < \infty.
\]
\end{proof}

So far, we have shown that $\bar{F}(i,j,t) \sim t^{-\alpha}L(t)
\left(\sum_{n \geq 0}n \Prob \{ N^{i,j}_0=n \} +
\delta_{i,0}\right)$ for $j \neq 0$. In view of Lemmas
\ref{domlemma} and \ref{lemma-convergence-Fbar}, the
representation (\ref{F*}) of $F^*(i,j,t)$ yields a similar result
for the stationary setting.

\begin{corollary} \label{cor-barF*} For all $i,j \neq 0$ we have
\[
     \lim_{t \rightarrow \infty} \frac{\bar{F}^*(i,j,t)}{t^{-\alpha}
     L(t)} = \sum_{n \geq 0}n \Prob \{ N^{i,j}_0=n \}.
\]
\end{corollary}

\begin{proof}
Due to (\ref{F*}), (\ref{defnofhats}) and (\ref{star}) we can
write
\begin{equation}\label{beatFstar}
F^{*}(i,j,t)=\frac{p_{ij}}{m_{i,j}}\int_{0}^{t}(1-G(i,j,s))ds+
\sum_{k \neq j} \frac{p_{ik}}{m_{i,k}} \int_{0}^{t} F(k,j,t-u)
(1-G(i,k,u))du,
\end{equation}
and therefore
\begin{equation}\label{barFstarintermsofG}
\begin{split}
\bar{F}^{*}(i,j,t)
 &=\frac{p_{ij}}{m_{i,j}} \int_{t}^{\infty}(1-G(i,j,s))ds+
\sum_{k \neq j} p_{ik} \left(1-\frac{1}{m_{i,k}}\int_{0}^{t}
F(k,j,t-u) (1-G(i,k,u))du\right).
\end{split}
\end{equation}
We will now show that
\begin{equation}\label{intoghatG}
\lim_{t \rightarrow \infty }
\frac{\int_{t}^{\infty}(1-G(i,j,s))ds}{t^{-\alpha}L(t)}=0 \quad
\mbox{if} \quad p_{ij} > 0 ~~ \mbox{ and } ~~ i \neq 0.
\end{equation}
To this end, we first apply Proposition 1.5.10 in \cite{bingham}:
if $g$ is a function on $\re_+$ that satisfies $g(t) \sim
t^{-\beta} L(t)$ for $\beta > 1$, then
\[
    \int_t^\infty g(s) ds \sim \frac{t^{-\beta + 1}}{\beta - 1}
    L(t).
\]
Together with Assumption~\ref{Assumption-tails} (iii) this
proposition implies that
\begin{equation}\label{eq:hdec}
\lim_{t \rightarrow \infty} \frac{\int_t^\infty
h(i,s)ds}{t^{-\alpha} L(t)}=0, \quad i \neq 0,
\end{equation}
where $h(i,\cdot)$ is the tail of the distribution of the sojourn
time in state $i$, defined in (\ref{hdef}). The representation
\begin{equation}\label{hrep}
h(i,t)= 1-\sum_{j \in E}p_{ij}G(i,j,t)= \sum_{j \in E} p_{ij}
(1-G(i,j,t)),
\end{equation}
along with (\ref{eq:hdec}) implies (\ref{intoghatG}).

In order to find the decay rate of the remaining terms of
(\ref{barFstarintermsofG}), recall first that $N^{i,j}_k$ is the
number of visits to state $k$ before reaching state $j$, given
$\xi_0 = i$, and not counting the first one in the case $k=i$.
\[
    \mathbb{P}\{N^{i,j}_0 = n\} = \sum_{k \notin\{0,j\}}p_{ik}
    \mathbb{P}\{N^{k,j}_0 = n\} + p_{i0} \mathbb{P}\{N^{i,0}_0 =
    n-1\},
\]
and so
\begin{equation} \label{A}
     \sum_{n \geq 0} n \mathbb{P}\{N^{i,j}_0 = n\} = \sum_{k \neq j}
    p_{ik} \sum_{n \geq 0} n \mathbb{P}\{N^{k,j}_0 = n\} + p_{i0} \sum_{n \geq 0}
    \mathbb{P}\{N^{0,j}_0=n\}.
\end{equation}
Note that Assumption \ref{Assumption-pi} implies $\sum_{n \geq 0}
    \mathbb{P}\{N^{0,j}_0=n\}=1$. Since
$\frac{1}{m_{i,j}}\int_{0}^{t}(1-G(i,j,s))ds$ is a distribution
function whose tail is
$\frac{1}{m_{i,j}}\int_{t}^{\infty}(1-G(i,j,s))ds$, Lemmas
\ref{domlemma}, \ref{lemma-convergence-Fbar} together with
equations (\ref{intoghatG}) and (\ref{A}) imply that
\begin{eqnarray*}
     \lim_{t \rightarrow \infty}\frac{\sum_{k \neq j} p_{ik}
    \left(1-\frac{1}{m_{i,k}}\int_{0}^{t} F(k,j,t-u)
    (1-G(i,k,u))du\right)}{t^{-\alpha}L(t)}
    &=& \sum_{k \neq j}p_{ik}
    \sum_{n\geq 0}n \mathbb{P}\{N_{0}^{k,j}=n\} + p_{i0}\\
    & =& \sum_{n \geq 0} n \mathbb{P}\{N_{0}^{i,j}=n\}.
\end{eqnarray*}
This completes the proof.

\end{proof}

The next result shows that the first term on the right-hand-side
of (\ref{eq:ipij}) converges to zero sufficiently fast.
\begin{corollary} \label{convergence-s}
For all $i \neq 0$ we have
\[
    \lim_{t \rightarrow \infty} \frac{s(i,t)}{t^{-\alpha + 1} L(t)} =
    0.
\]
\end{corollary}

\begin{proof}
The proof is an immediate consequence of Corollary~\ref{cor-barF*}
since $\bar{F}^{*}(i,j,t) \geq s(i,t)/\nu_i$.
\end{proof}



\subsubsection{The tail structure of $R^* * h$}

So far, we have analyzed the tail structure of the distribution of
the travel time between states $i$ and $j$ $(i,j \neq 0)$ in the
stationary regime. We are now going to study the tail structure of
$R^*(i,j,t) * h(j,t)$.

\begin{lemma} \label{convergence-R*h}
There exists $C_j > 0$ such that, for $j \neq 0$,
\[
    \lim_{t \rightarrow \infty} \frac{R^*(i,j,t)*h(j,t) - \nu_j}
    {t^{-\alpha + 1} L(t)} = \frac{C_j}{\alpha - 1}
\]
for all $i \in E$, $i \neq 0$.
\end{lemma}
\begin{proof}
Let us fix $i,j \in E$, $i,j \neq 0$. Using the representation
(\ref{eq:srij}) for the function $R^*(i,j,\cdot)$, we need to show
that
\[
    \lim_{t \rightarrow \infty} \frac{R(j,j,t)*F^*(i,j,t)* h(j,t) -
    \nu_j} {t^{-\alpha + 1} L(t)} = \frac{C_j}{\alpha-1}.
\]
It follows from Fubini's theorem that
\begin{eqnarray*}
    \int_{0}^{\infty} F^{*}(i,j,t)* h(j,t) dt & = &
    \int_{0}^{\infty} \int_{0}^{\infty} h(j,t-s)\textbf{1}_{\{s \leq t\}} F^*(i,j,ds)\, dt \\
    & = & \int_{0}^{\infty} \int_{s}^{\infty} h(j,t-s) dt\, F^*(i,j,ds) \\
    & = & \int_{0}^{\infty} h(j,t)dt \\
    & = & m_j,
\end{eqnarray*}
where $m_j$ is the mean sojourn time at state $j \in E$ as defined
in (\ref{mean-time}). 
Suppose now that we can show that the continuous non-negative
function $z$ defined by
\begin{equation} \label{def-z}
    z(i,j,t) := F^*(i,j,t)*h(j,t)
\end{equation}
satisfies $z(i,j,t) = o(\bar{F}(i,j,t))$, is directly Riemann
integrable\footnote{See Appendix \ref{Appendix-Riemann} for the
definition of directly Riemann integrability.} and of bounded
variation, and its total variation function $z^{*}(i,j,\cdot)$,
defined in (\ref{z-star}) below, satisfies
$z^{*}(i,j,t)=O(t^{-\alpha+1}\hat{L}(t))$, where we define
\[ \hat{L}(t) =L(t) \sum_{n \geq 0} n \Prob\{N^{j,j}_0 = n \}. \]
Then Theorem \ref{sidthm} yields
\begin{eqnarray} \label{tail-z}
    \frac{m_j}{\eta_j} - \int_0^t z(i,j,t-s) R(j,j,ds)
    & \sim & - \frac{m_j}{(\alpha-1)\eta_j^2} t^{-\alpha + 1}
    \hat{L}(t) ,
\end{eqnarray}
because $R(j,j,t) = \sum_{n \geq 0} F^{n}(j,j,t)$ is a renewal
function, $F(j,j,\cdot)$ is nonsingular and $ \bar{F}(j,j,t) \sim
t^{-\alpha} \hat{L}(t)$ for $j \neq 0$.
By Propositions 5.5 and 6.12 in \cite{Cinlar75}, $\nu_j =
\frac{m_j}{\eta_j}$, and so (\ref{tail-z}) would yield
\[
    \lim_{t \rightarrow \infty} \frac{R(j,j,t)*F^{*}(i,j,t)*h(j,t) - \nu_j}{t^{- \alpha
    + 1} L(t)} =  \frac{m_j}{(\alpha - 1)\eta_j^2}
    \sum_{n \geq 0} n \Prob\{N^{j,j}_0 = n \}
\]
and hence prove our assertion with
\begin{equation}\label{Cj}
    C_j:=  \frac{m_j}{\eta^2_j}
    \sum_{n \geq 0} n \Prob\{N^{j,j}_0 = n \}.
\end{equation}
We will now establish that $(i)~z(i,j,t) = o(\bar{F}(i,j,t))$, that
$(ii)~z^{*}(i,j,t)=O(t^{-\alpha+1}\hat{L}(t))$ and $z(i,j,\cdot)$ is
of bounded variation, and that $(iii)~z(i,j,\cdot)$ is
directly-Riemann integrable.

\begin{itemize}
\item[(i)]
Since $i,j \neq 0$ and because $p_{i,0} > 0$ the probability that
the semi-Markov process $x$ visits the state $0$ before it reaches
the state $j$ is positive. Thus, it follows from Assumption
\ref{Assumption-tails} (iii) and from Corollary \ref{cor-barF*}
that
\begin{equation}\label{decayofh}
    h(j,t) = o(\bar{F}^*(i,j,t)).
\end{equation}
Now it follows from (\ref{decayofh}) and Lemma \ref{domlemma} that
\begin{equation} \label{bounded-variation}
\bar{F}^*(i,j,t)+z(i,j,t)=1-F^{*}(i,j,t)*(1-h(j,t)) \sim
\bar{F}^{*}(i,j,t).
\end{equation}
Hence $z(i,j,t)=o(\bar{F}^{*}(i,j,t))$, and so
$z(i,j,t)=o(\bar{F}(i,j,t))$ by Lemma \ref{lemma-convergence-Fbar}
and Corollary \ref{cor-barF*}.

\item[(ii)] Assumption~\ref{Assumption-tails} (iv) along with
(\ref{eq:F2}), the representation (\ref{beatFstar}) of
$F^*(i,j,\cdot)$ and (\ref{hrep}) guarantees that $F^*(i,j,\cdot)$
and $h(j, \cdot )$ have a bounded continuous densities
$f^{*}(i,j,\cdot)$ and $h'(j, \cdot)$, respectively. As a result,
$z(i,j,\cdot)$ is absolutely continuous with density
\begin{equation} \label{density-z}
    z'(i,j,t)= f^*(i,j,t) + \int_0^t
    h'(j,t-s)f^*(i,j,s)ds.
\end{equation}
Since $h'(j,t-s) \leq 0$,
\begin{equation} \label{abs-z}
    |z'(i,j,t)| \leq f^*(i,j,t)-\int_{0}^{t}
    h'(j,t-s) f^{*}(i,j,s) ds.
\end{equation}
In terms of the distribution function $g(j,t) := 1 - h(j,t)$ we
have
\begin{equation}\label{lastmohk}
\begin{split}
    - \int_{0}^{t} h'(j,t-s) f^{*}(i,j,s)ds
    & = - \frac{\partial}{\partial t} \int_{0}^{t} h(j,t-s)
    f^{*}(i,j,s) ds + f^*(i,j,t) \\
    & = \frac{\partial}{\partial t}
    \left(g(j,t)*F^{*}(i,j,t)\right).
\end{split}
\end{equation}
As a result the total variation function $z^*(i,j,\cdot)$
satisfies
\begin{equation} \label{z-star}
    z^*(i,j,t) = \int_t^\infty |z'(i,j,s)| ds
    \leq \bar{F}^*(i,j,t) + \left[ 1 - g(j,t) * F^*(i,j,t)
    \right].
\end{equation}
Thus, $z^*(i,j,t) \leq 2$ and $z(i,j,\cdot)$ is of bounded
variation. Lemma~\ref{domlemma} and (\ref{decayofh}) together with
(\ref{z-star}) give $z^{*}(i,j,t) = O(\bar{F}^{*}(i,j,t))$, and so
Corollary~\ref{cor-barF*} yields
$z^{*}(i,j,t)=O(t^{-\alpha+1}\hat{L}(t))$.

\item[(iii)]
In order to prove that $z$ is directly Riemann integrable note
first that since $h \geq 0$ is decreasing and by
Assumption~\ref{Assumption-tails} (i) it is integrable, therefore
by Lemma~\ref{cinlarslemma} (i) it is directly Riemann integrable.
Now Lemma~\ref{cinlarslemma} (ii) implies that $z$ is directly
Riemann integrable.

\end{itemize}
\end{proof}


\subsubsection{Proof of the central limit theorem for semi-Markov
processes}

We are now ready to prove the main result of this section.

\mbox{ }

\noindent \textsc{Proof of Theorem \ref{thm2}:} By (\ref{eq:ipij})
we have the representation
\[
    P^*_{t}(i,j)=\frac{s(i,t)}{\nu_{i}}\delta_{ij}+\int_{0}^{t}
    h(j,t-s) R^*(i,j,ds)
\]
for the conditional probability that $x_t=j$, given $x_0=i$. Due
to Corollary \ref{convergence-s} and Lemma \ref{convergence-R*h},
\[
    \lim_{t \rightarrow \infty} \frac{P^*_{t}(i,j) - \nu_j}{t^{-\alpha
    +1}L(t)} = \frac{C_j}{\alpha - 1} \qquad (i,j \neq 0),
\]
where $C_{j}$ is given by (\ref{Cj}). With $H=\frac{3-\alpha}{2}$
it follows from (\ref{eq:cov}) that
\begin{equation}\label{gammin}
    \lim_{t \rightarrow \infty} \frac{\gamma(t)}{t^{2H-2} L(t)} =
    \frac{1}{(2-2H)} \sum_{i,j \in E} i j \nu_i C_j.
\end{equation}
By Lemma \ref{lemma-covariance} this proves the existence of a
constant $c$ such that the finite dimensional distributions of the
processes $\left(\frac{1}{\eps^{1-H} \sqrt{L(\eps^{-1})}}
\,X^{\eps}_{t} \right)_{0 \leq t < \infty}$ converge weakly to the
finite dimensional distributions of the fractional Brownian motion
process $cB^H$ as $\eps\downarrow 0$, and $c$ is given by
\[
    c^2 = \frac{1}{2H(1-H)(2H-1)}\sum_{i,j \in
    E}ij\nu_{i}\frac{m_{j}}{\eta_{j}^{2}}  \sum_{n \geq 0} n \Prob\{N^{j,j}_0 = n \} = \frac{1}{2H(1-H)(2H-1)}
    \mu \sum_{j \in E} j \frac{m_{j}}{\eta_{j}^{2}}  \sum_{n \geq 0} n \Prob\{N^{j,j}_0 = n \}.
\]

In order to establish tightness, we proceed in two steps.

\begin{rmenumerate}
\item We first establish the existence of a constant $C <
\infty$ such that
\begin{equation} \label{all-t}
    \Var(t) \leq C t^{2H} L(t) \qquad \mbox{for all} \qquad t \geq
    0.
\end{equation}
In view of (\ref{variance}), we can choose a sufficiently large $T
\in \n$ that satisfies
\[
    \Var(t) \leq 2 c^2 t^{2H}L(t) \qquad \mbox{for all} \qquad t
    \geq T.
\]
In terms of the random variables $Y^i := \int_{(i-1)t/T}^{it/T}
Y_s ds$ $(i=1,2, \ldots ,T)$ we have
\[
    \Var\left( \int_0^t Y_s ds \right) = \Var\left( \sum_{i=1}^T
    \int_{(i-1)t/T}^{it/T} Y_s ds \right) = \sum_{i,j=1}^T
    \textnormal{Cov}(Y^i,Y^j).
\]
Since the mean zero Gaussian process $(\int_0^\cdot Y_s ds)$ has
stationary increments, the sequence $Y^1, Y^2, \ldots ,Y^T$ is
stationary. Thus, the H\"older inequality yields
\[
    \textnormal{Cov}(Y^i,Y^j) \leq \sqrt{\Var(Y^i)}
    \sqrt{\Var(Y^j)} = \Var(Y^0),
\]
and so
\[
    \Var(t) \leq T^2 \Var\left( \frac{t}{T} \right) \quad
    \mbox{for} \quad 0 \leq t \leq T.
\]

In view of (\ref{eq:cov}), the function $\gamma$ is bounded:
$\|\gamma\|_\infty < \infty$. Hence for $s \in [0,1]$ the
representation (\ref{representation-variance}) shows that
\[
    \Var(s) = 2 \int_0^s \int_0^v \gamma(u) du dv \leq
    s^2 \|\gamma\|_\infty \leq s^{2H} \|\gamma\|_\infty \quad
    \mbox{because} \quad H \in \left( \frac{1}{2}, 1 \right).
\]
This yields (\ref{all-t}) with $C := \max\{2 c^2, T^{2-2H}
\|\gamma\|_\infty\} $ (after putting $L \equiv 1$ on $[0,T]$).

\item Let us now denote by $Z^{\eps}=(Z^{\eps}_t)$ the mean zero
Gaussian process with stationary increments defined by
\[
    Z^{\eps}_t :=
    \frac{1}{\eps^{1-H}\sqrt{L(\eps^{-1})}}
    X^{\eps}_t \qquad (0 \leq t \leq T).
\]
Due to Theorem 12.3 in \cite{Billingsley}, the family of
stochastic processes $(Z^{\eps})$ is tight if the following moment
condition is satisfied for some constants $\delta > 1$ and
$\hat{C} < \infty$ and for all sufficiently small $\eps$:
\begin{equation}\label{tightness}
    \E[ Z^\eps_{t_2} - Z^\eps_{t_1}]^2
    \leq \hat{C}|t_{2}-t_{1}|^{\delta} \quad \mbox{for all} \quad 0 \leq t_1
    \leq t_2 \leq T.
\end{equation}
In view of (\ref{Xepsdef}), we have
\[
    \E [Z^\eps_{t_{2}} - Z^\eps_{t_{1}} ]^{2} =
    \eps^2 \Var\left(\frac{t_{2}-t_{1}}{\eps}\right)
    \frac{1}{\eps^{2-2H}L(\eps^{-1})},
\]
and so it follows from step (i) that
\begin{equation} \label{estimate-increment}
   \E [Z^\eps_{t_{2}} - Z^{\eps}_{t_{1}} ]^{2} \leq
   C (t_{2}-t_{1})^{2H} \frac{L\left(\eps^{-1}(t_{2}-t_{1}) \right)}
   {L(\eps^{-1})}.
\end{equation}
Since $L$ is slowly varying, $\frac{L(\eps^{-1}u)}{L(\eps^{-1})}$
tends to 1 as $\eps \rightarrow 0$, and this convergence is
uniform in $u$ over compact sets (\cite{bingham}, Theorem 1.2.1).
Thus, there exists $\eps^*>0$ such that
\[
    \frac{L(\eps^{-1}u)}{L(\eps^{-1})} \leq 2 \quad \mbox{for all}
    \quad 0 \leq u \leq T
\]
if $\eps < \eps^*$, and so (\ref{tightness}) follows from
(\ref{estimate-increment}) with $\hat{C} := 2 C$. \hfill $\Box$
\end{rmenumerate}

\section{Approximating Integrals with respect to
\fbm}\label{Section-sigma}

In this section we prove an approximation result for stochastic
integrals which contains Theorem \ref{thm1} as a special case.
More precisely we give conditions which guarantee that for a
sequence of processes $\{(\Psi^n,Z^n)\}_{n \in \n}$ the
convergence $\mathcal{L}\mbox{-}\lim_{n \rightarrow \infty}
(\Psi^n,Z^n) = (\Psi,B^H)$ implies the convergence
$\mathcal{L}\mbox{-}\lim_{n \rightarrow \infty} \left( \Psi^n,Z^n,
\int \Psi^n dZ^n \right) = \left( \Psi,B^H, \int \Psi dB^H
\right)$. We follow the notational convention in
\cite{Kurtz-Protter} that $\int X dY$ (without displaying the
running variable) denotes $\int X_{s-}dY_s$.

All stochastic integrals in this section are understood as
 limits in probability of Stieltjes-type sums which can be described as follows: given a
 probability space $(\Omega, \mathcal{F}, \mathbb{P})$ and a
 filtration $\{\mathcal{F}_{t}: t \geq 0\}$ of sub-sigma-fields of $\mathcal{F}$
consider two adapted stochastic processes $\phi$ and $Z$; we say
that the integral $\int \phi dZ$ exists if for any $T < \infty$
and for each sequence of partitions $\{\tau^l\}_{l \in \n}$,
$\tau^l = (\tau^l_1, \tau^l_2, \ldots ,\tau^l_{N_l})$, of the
interval $[0,T]$ that satisfies $\lim_{l \rightarrow \infty}
\max_i |\tau^l_{i+1} - \tau^l_i | = 0$,
\begin{equation}
\label{definition-integral}
    \int_0^T \phi_{s-} dZ_{s} = \Prob\mbox{-}\lim_{l \rightarrow \infty}
    \sum_i \phi_{\tau^l_i}(Z_{\tau^l_{i+1}} - Z_{\tau^l_{i}}),
\end{equation}
where $\Prob\mbox{-}\lim$ denotes the limit in probability.
 This definition of stochastic integrals applies to the
usual semimartingale setting where $\phi$ is a process in
$\mathbb{D}$ and where $Z$ is a semimartingale. If $Z = B^H$ is a
fractional Brownian motion process with Hurst coefficient $H >
\frac{1}{2}$, the limit in (\ref{definition-integral}) exists for
a large class of integrands, including continuous semimartingales
and $C^{1}$-functions of fractional Brownian motion. In
particular, the stochastic integral $\int B^H dB^H$ exists in the
sense of (\ref{definition-integral}), and for the continuous
semimartingale $\Psi$ we have the following integration by parts
formula, due to \cite{Lin}:
\begin{equation}
\label{integration-by-parts}
    - \int B^{H} d \Psi + \Psi B^H = \int \Psi d B^H.
\end{equation}

Before we state the main result of this section, we recall that a
sequence $\{\Psi^n\}_{n \in \n}$ of semimartingales defined on
probability spaces $(\Omega^n, \mf^n, \Prob^n)$ is called
\textsl{good} in the sense of \cite{Duffie-Protter} if, for any
sequence $\{Z^n\}_{n \in \n}$ of c$\grave{a}$dl$\grave{a}$g
adapted processes, the convergence $\mathcal{L}\mbox{-}\lim_{n
\rightarrow \infty} (\Psi^n,Z^n) = (\Psi,Z) $ implies the
convergence
\[
    \mathcal{L}\mbox{-}\lim_{n \rightarrow \infty}
    \left(\Psi^n,Z^n,\int Z^n d\Psi^n \right)
    = \left(\Psi,Z,\int Z d \Psi \right).
\]

\begin{remark} \label{dufprolem}
Any semi-martingale $\Psi$ can be written as $\Psi=M+A$, in which
$M$ is a local martingale with $M_{0}=0$ and $A$ is an adapted
finite variation process. We will denote by $[M,M]$ the quadratic
variation of $M$ and by $|A|$ the total variation of $A$. %
%
A sequence $(\Psi^{n}\}_{n \in \n}$ $(\Psi^n = M^{n}+A^n)$ of
semi-martingales is good on $[0,T]$ if the sequences
$(\mathbb{E}_{n}\{[M^{n},M^{n}]_{T}\})_{n \in \mathbb{N}}$ and
$(\mathbb{E}_{n}\{|A|_{T}\})_{n \in \mathbb{N}}$ are bounded; see
\cite[Theorem 4.1]{Duffie-Protter} for details.
\end{remark}

In the following we will denote by $\{\Psi^n\}_{n \in \n}$ a
sequence of ``good'' semimartingales and by $\{Z^n\}_{n \in \n}$ a
sequence of $\mathbb{D}$-valued stochastic processes defined on
some probability space $(\Omega,\mf,\Prob)$ and assume that the
following conditions are satisfied.

\begin{assumption} \label{Assumption3}
\begin{rmenumerate}
\item The sample paths of the processes $Z^n$ are almost surely of
zero quadratic variation on compact sets, and $\Prob \{Z^n_0 = 0\}
= 1$. \item The stochastic integrals $\int \Psi^n d Z^n$ and $\int
Z^n dZ^n$ exist in the sense of (\ref{definition-integral}), and
the sample paths $t \mapsto \int_0^t Z^n_{s-} d Z^n_s$ and $t
\mapsto \int_0^t \Psi^n_{s-} d Z^n_s$ are c\`{a}dl\`{a}g.
\end{rmenumerate}
\end{assumption}

We are now ready to state the main theorem of this section. Its
proof requires some preparation which will be carried out below.

\begin{theorem} \label{thm3} Let $\{\Psi^n\}_{n \in \n}$ be a
sequence of good semimartingales and let $\{Z^n\}_{n \in \n}$ be a
sequence of $\mathbb{D}$-valued stochastic processes that satisfy
Assumption \ref{Assumption3}. If $\Psi$ is a continuous
semimartingale and if $B^H$ is a \fbm process with Hurst parameter
$H > \frac{1}{2}$, then the convergence
$\mathcal{L}\mbox{-}\lim_{n \rightarrow \infty}(\Psi^n,Z^n) =
(\Psi,B^H)$ implies the convergence
\[
    \mathcal{L}\mbox{-}\lim_{n \rightarrow \infty}
    \left( \Psi^n,Z^n,\int \Psi^n dZ^n \right) =
    \left( \Psi,B^H, \int \Psi d B^H \right).
\]
\end{theorem}

\mbox{ }

Before we turn to the proof of Theorem \ref{thm3}, we consider an
example where Assumption \ref{Assumption3} can indeed be verified.

\begin{example} \label{example1}
    Let $\{H_n\}_{n \in \n}$ be a sequence of real numbers with $H_n >
    \frac{1}{2}$, and assume that $\lim_{n \rightarrow \infty}H_n = H >
    \frac{1}{2}$. Let $Z^n$ be a fractional Brownian motion process with Hurst
    parameter $H_n$ and let $\Psi$ be a continuous semimartingale independent of $Z^{n}$ for all $n$.
    Since $H_n > \frac{1}{2}$, the processes $Z^n$ have zero quadratic variation.
    Moreover, $\mathcal{L}\mbox{-}\lim_{n \rightarrow \infty}
    Z^n = B^H$ because the centered Gaussian processes $Z^n$ and $B^H$ are uniquely determined
    by their covariation functions and
    all stochastic integrals exists in the sense of (\ref{definition-integral}). Thus,
    Theorem \ref{thm3} yields
\[
    \mathcal{L}\mbox{-}\lim_{n \rightarrow \infty} \left( Z^n, \int \Psi d Z^n
    \right) = \left( B^H, \int \Psi d B^H \right).
\]
\end{example}

We prepare the proof of Theorem \ref{thm3} with the following
simple lemma.

\begin{lemma} \label{lemma-covariance2}
Under the assumptions of Theorem \ref{thm3} the processes $[Z^n]$
and $[Z^n,\Psi^n]$ defined by
\[
    [Z^n]_t := (Z^n_t)^2 - 2 \int_{0}^{t} Z^n_{s-} dZ^n_{s} \quad
    \mbox{and} \quad
    [Z^n,\Psi^n]_t := Z^n_{t} \Psi^n_t - \int_0^t Z^n_{s-} d \Psi^n_s
    - \int_0^t \Psi^n_{s-} d Z^n_s,
\]
have $\Prob$-a.s.~sample paths which are equal to zero.
\end{lemma}
\begin{proof}
It follow from the representation of the stochastic integrals
$\int Z^n dZ^n$, $\int \Psi^n dZ^n$ and $\int Z^n d \Psi^n$ as
probabilistic limits of Stieltjes-type sums that, for any $t$ and
each sequence of partitions $\{\tau^l\}_{l \in \n}$ of $[0,t]$
with $\lim_{l \rightarrow \infty} \max_i | \tau^l_{i+1} -
\tau^l_i|$,
\begin{equation*}
    [Z^n] = \Prob\mbox{-}\lim_{l \rightarrow \infty}
    \sum_i (Z^n_{\tau^l_{i+1}} - Z^n_{\tau^l_{i}})^2 ~~ \mbox{and} ~~
    [Z^n,\Psi^n] = \Prob\mbox{-}\lim_{l \rightarrow \infty}
    \sum_i (Z^n_{\tau^l_{i+1}} - Z^n_{\tau^l_{i}})
    (\Psi^n_{\tau^l_{i+1}} - \Psi^n_{\tau^l_{i}}).
\end{equation*}
Since a typical sample path of the stochastic integrals $\int
\Psi^n dZ^n$ and $\int Z^n d \Psi^n$ is in $\mathbb{D}$, we can
apply the same arguments as in the proof of Theorem II.6.25 in
\cite{Protter} in order to obtain the inequality $\Prob \left\{
[Z^n,\Psi^n]_t^2 \leq [Z^n]_t [\Psi^n]_t \right\} = 1$. Thus, our
assertion follows from $\Prob\{[Z^n]_t = 0 \mbox{ for all } t \geq
0 \} = 1$.
\end{proof}

For the proof of Theorem \ref{thm3} we will also need the
following result.

\begin{lemma} \label{lemma-D}
\begin{rmenumerate}
\item Let $\mathbb{C}$ be the space of all real valued continuous
functions. For $n \in \n$, let $\alpha_n, \beta_n \in \mathbb{D}$
and assume that the sequence $\{(\alpha_n,\beta_n)\}_{n \in \n}$
converges in the Skorohod topology to $(\alpha,\beta) \in
\mathbb{C}\times\mathbb{C}$. Then, on compact intervals, the
process
\[
    \gamma_n = (\gamma_n(t))_{t \geq 0}
    \quad \mbox{defined by} \quad \gamma_n(t) =
    \alpha_n(t) \beta_n(t)
\]
converges to $\alpha \beta = (\alpha(t) \beta(t))_{t \geq 0}$ in
the Skorohod topology on $\mathbb{D}$.
\item  Let $\{(Y^n,Z^n)\}_{n \in \n}$ be a sequence of
$\mathbb{D}$-valued random variables defined on some probability
space $(\Omega,\mf,\Prob)$ that converges in law to $(Y,Z)$. If
$\Prob\{(Y,Z) \in \mathbb{C} \times \mathbb{C}\}  = 1$, then
\[
    \mathcal{L}\mbox{-}\lim_{n \rightarrow \infty} \{(Y^n_t Z^n_t)_{0
    \leq t \leq T}\} = (Y_t Z_t)_{0 \leq t \leq T}
\]
holds for all $T < \infty$.
\end{rmenumerate}
\end{lemma}
\begin{proof} Since $\alpha$ and $\beta$ are continuous, (i)
follows from Lemma 2.1 in \cite{Kurtz-Protter}. The second
assertion follows from (i) and Skorohod's representation theorem.
\end{proof}

We are now ready to finish the proof of Theorem \ref{thm3}.

\mbox{ }

\noindent{\textsc{Proof of Theorem \ref{thm3}}:} Since
$\{\Psi^n\}_{n \in \n}$ is a sequence of good semimartingales and
because a typical sample path of a fractional Brownian motion
process is continuous, we deduce from Theorem 2.2 in
\cite{Kurtz-Protter} and from Lemma \ref{lemma-D} (ii) that
\begin{equation}
\label{limits}
    \mathcal{L}\mbox{-}\lim_{n \rightarrow \infty} \left( \Psi^n, Z^n,
    \int Z^n d\Psi^n \right) = \left( \Psi, B^H, \int B^H d \Psi \right)
    ~ \mbox{and} ~
    \mathcal{L}\mbox{-}\lim_{n \rightarrow \infty} \left( \Psi^n Z^n
    \right) = \left( \Psi B^H \right),
\end{equation}
respectively. By the continuous mapping theorem, it follows from
(\ref{limits}), from Lemma \ref{lemma-covariance2} and from the
integration by parts formula for fractional Brownian motion
(\ref{integration-by-parts}) that the finite dimensional
distributions of the processes
\[
    \left( \Psi^n, Z^n, \int \Psi^n d Z^n \right) =
    \left( \Psi^n, Z^n, - \int Z^n d \Psi^n + \Psi^n Z^n
    \right)
\]
converge weakly to the finite dimensional distributions of the
process
\[
    \left( \Psi, B^H, - \int B^H d \Psi + \Psi B^H \right) =
    \left( \Psi, B^H, \int \Psi d B^H \right).
\]
It also follows from (\ref{limits}) that the sequence $\left\{
\int Z^n d\Psi^n \right\}_{n \in \n}$ is C-tight and the sequence
$\left\{ \Psi^n Z^n \right\}_{n \in \n}$ is tight. By Corollary
VI.3.33 in \cite{JS} the sum of a tight sequence of stochastic
processes with a sequence of C-tight processes is tight. Thus,
continuity of the processes $\Psi B^H$ and $\int B^H d \Psi$
yields tightness of the sequence $\left\{ \int \Psi^n dZ^n
\right\}_{n \in \n}$. This shows that
\[
    \mathcal{L}\mbox{-}\lim_{n \rightarrow \infty} \left( \Psi^n,
    Z^n, \int \Psi^n d Z^n \right) = \left( \Psi, B^H, \int \Psi d B^H
    \right).
\]
\hfill $\Box$

\medskip

In view of Lemma \ref{lemma-covariance2} and the integration by
parts formula for fractional Brownian motion processes we also
have the following approximation result for the integral of
fractional Brownian motion process with respect to itself.

\begin{proposition}
    Under the assumption of Theorem \ref{thm3} it holds that
\[
    \mathcal{L}\mbox{-}\lim_{n \rightarrow \infty} \int_{0}^{t} Z^n_{s-}
    dZ^n_s
    = \int_{0}^{t} B^H_s d B^H_s.
\]
\end{proposition}
\begin{proof} By Lemma \ref{lemma-covariance2}
\[
    (Z^n_t)^2 = 2 \int_0^t Z^n_{s-} d Z^n_s \qquad \Prob\mbox{-a.s.}
\]
Thus, in view of Lemma \ref{lemma-D} (ii) and the It\^{o} formula
for fractional Brownian motion the sequence $\{(Z^n_t)^2\}_{n \in
\n}$ converges in distribution to
\[
    (B^H_t)^2 = 2 \int_0^t B^H_s dB^H_s.
\]
This yields the assertion.
\end{proof}

We finish this section with the proof of Theorem \ref{thm1}.

\medskip

\noindent \textsc{Proof of Theorem \ref{thm1}:} In terms of the
processes $Y$ and $X^\eps$ introduced in (\ref{def-Y}) and
(\ref{Xepsdef}), respectively, we have
\begin{eqnarray*}
    \mathcal{L}\mbox{-}\lim_{N \rightarrow \infty}
    \frac{1}{\sqrt{N}}(X^{\eps,N}_t)_{0 \leq t \leq T} & = &
    \mathcal{L}\mbox{-}\lim_{N \rightarrow \infty}\left(
    \int_0^t \left(\frac{1}{\sqrt{N}} \sum_{a \in \ma} \Psi_s
    x^a_{s/\eps} - \sqrt{N} \mu \Psi_s\right)ds \right)_{0 \leq t \leq T} \\
    & = & \left( \int_0^t \Psi_s Y_{s/\eps} ds
    \right)_{0 \leq t \leq T} \\
    & = &
    \left( \int_0^t \Psi_s d X^\eps_s \right)_{0 \leq t \leq T}
\end{eqnarray*}
and
\[
    \int_0^t X^\eps_s dX^\eps_s =  \int_0^t X^\eps_s
    Y_{s/\eps}\,ds
\]
where all the stochastic integrals have continuous sample paths.
By Theorem \ref{thm2},
\[
    \mathcal{L}\mbox{-}\lim_{\eps\downarrow 0} \frac{1}{\eps^{1-H}\sqrt{L(\eps^{-1})}} X^\eps = c\,B^H
\]
with $H > \frac{1}{2}$ and the sequence of semi-martingales
consisting of a single element $\Psi$ is good as a result of
Remark~\ref{dufprolem} since by assumption
$E\{[\Psi,\Psi]_{T}\}<\infty$ and $\mathbb{E}\{|A|_{T}\}<\infty$.
Therefore, the assertion follows from Theorem \ref{thm3} if we can
show that the processes $X^\eps$ have zero quadratic variation on
compact time intervals and that the stochastic integrals $\int
\Psi d X^\eps$ and $\int X^\eps\, dX^\eps$ exist as the
probabilistic limits of Stieltjes-type sums. These properties,
however, follow from (\ref{Xepsdef}) by direct computation. \hfill
$\Box$


\section{Conclusion \& Discussion\label{conc}}
We proved a functional central limit theorem for stationary
semi-Markov processes and an approximation result for stochastic
integrals of fractional Brownian motion. Our motivation was to
provide a mathematical framework for analyzing microstructure
models where stock prices are driven by the demand of many small
investors. We identified inertia as a possible source of long
range dependencies in financial time series and showed how
investor inertia can lead to stock price models driven by a
fractional Brownian motion. Several avenues are open for further
research.

In our financial market model, investors do not react to changes
in the stock price. Over short time periods, such an assumption
might be justified for small, non-professional investors, but
incorporating feedback effects, whereby traders' investment
decisions can be influenced by asset prices, would be an important
next step from an economic point of view. This, however, leads to
a significant increase in the complexity of the dynamics, as
discussed in Remark \ref{twofour}.



Financial market models where stock prices are driven by the
demand of many heterogenous agents constitute another mathematical
challenge. Heterogeneity among traders has been identified as a
key component affecting the dynamics of stock prices, and in
recent years an extensive literature on the mathematical modelling
of heterogeneity and interactions in financial markets has
appeared. But the analysis is usually confined to models of active
market participants trading every period, and it is therefore of
interest to introduce heterogeneity and interaction into our
model. One approach would be to build on the case study in Section
\ref{twothree} to look at models where inert traders (who alone
would lead to a limit fractional Brownian motion price process)
interact with active traders (who alone would lead to a limit
standard Brownian motion price process), and to include mechanisms
by which their {\em interaction} leads to efficient markets. One
could also try incorporate {\em strategically} interacting
institutional investors; see \cite{bp} for a possible game
theoretic framework.





For long-term economic models, it is important to look at
time-varying rates of long-range dependence, where the variation
is caused by global economic factors, or regime changes. At some
times, the market may be efficient, for example in a bullish
exuberant economy like the late 1990s (see Figure \ref{Hests}),
while at other times, the Joseph effect may be prominent, such as
during a recession or period of economic nervousness as in the
early 1990s. The mathematical challenge is then to adapt the limit
theorem of this paper to stochastically varying measures of
inertia.



\begin{appendix}
\section{The key renewal theorem in the heavy tailed case}
\label{Appendix-renewal}

In this appendix we recall a result of Heath et al. (\cite{heath})
on the rate of convergence in the key renewal theorem in the heavy
tailed case.

\begin{theorem}\label{sidthm}
Let $F$ be a distribution with domain $[0,\infty)$ satisfying
\[
    \bar{F}(t) = 1-F(t) \sim t^{-\alpha} \hat{L}(t)
\]
for some $1 < \alpha <2$, and where $\hat{L}$ is a slowly varying
function at infinity. Assume that $F^{n}$ is nonsingular for some
$n \geq$ 1. Let $\kappa = \int_{0}^{\infty}\bar{F}(x)dx$ be the
expected value and denote by $U$ the renewal function associated
with $F$, that is,
\[
    U=\sum_{n=0}^{\infty}F^{n}.
\]
Let $z$ be a continuous, non-negative function of bounded
variation on $[0,\infty)$, such that $\lim_{t \rightarrow \infty}
z(t) = 0$. That is, $z(t) = \int_t^\infty \zeta(dy)$ for some
finite signed measure $\zeta$ on $[0,\infty)$. Let $z^*$ denote
the total variation function of $\zeta$. That is,
\begin{equation*}
    z^*(t) = \int_t^\infty |\zeta| (dy).
\end{equation*}
We will also assume that $z$ is a \emph{directly Riemann
integrable function} (see \cite{cinlar}, p.295 for a definition)
on $[0,\infty)$, such that $z(t) = o(\bar{F}(t))$
and that
\begin{equation}\label{zscond}
z^{*}(t)=O \left( t^{-\alpha+1} \hat{L}(t) \right). 
\end{equation}
Let $\lambda =\int_{0}^{\infty}z(t)dt < \infty$. Then the function
$h: \re_+ \rightarrow \re_+$ defined by
\[
    h(t) = \frac{\lambda}{\kappa}-\int_{0}^{t}z(t-s)U(ds)
\]
satisfies
\[
    h(t) \sim - \frac{\lambda}{(\alpha-1)\kappa^2} t^{-\alpha +
    1} \hat{L}(t).
\]
\end{theorem}
\begin{proof}
The proof follows from the Remark on page 11 of \cite{heath}.
\end{proof}

\section{Directly Riemann Integrable Functions} \label{Appendix-Riemann}

The proof of our approximation result uses the notion of directly
Riemann integrability. A real valued function $g$ on $\re_+$ is
called directly Riemann integrable if
\[
    \lim_{a \rightarrow 0+} \sum_{n=1}^\infty a \inf_{(n-1)a \leq t
    \leq na} g(t)
    = \lim_{a \rightarrow 0+}
    \sum_{n=1}^\infty a \sup_{(n-1)a \leq t \leq na} g(t)
    = \int g(t) dt.
\]
These sums would be Riemann sums if not for the infinite limit of
summation. If a function is directly Riemann integrable, it is
Lebesgue integrable, but the converse is not necessarily true.

\begin{lemma}\label{cinlarslemma}(\cite[Chapter 9]{cinlar}, Proposition 2.16, (c) and (d))
\begin{itemize}
\item[(i)] Let $g \geq 0$ be a monotone non-increasing function;
then $g$ is directly Riemann integrable if and only if g is
Riemann integrable.
\item[(ii)] Let $g \geq 0$ and let $\phi$ be a distribution
function on $\mathbb{R}_{+}$. If $g$ is directly Riemann
integrable, then $\phi*g$ is directly Riemann integrable as well.
\end{itemize}

\end{lemma}

\section{A Goodness Criterion for Semi-Martingales}

%

\end{appendix}

\acknowledgments{ We would like to thank the anonymous referees
for their insightful comments, which helped to improve the paper.
The work of the first author supported by the US Office of Naval
Research under grant no. N00014-03-1-0102. Work of the second
author supported by the German Academic Exchange Service, the DFG
research center ``Mathematics for key technologies'' (FZT 86) and
by NSERC grant no. 312600-05. We would like to thank Peter Bank,
Erhan \c{C}inlar, Hans F\"ollmer, and Chris Rogers for valuable
comments. Part of this work was carried out while the second
author was visiting the Bendheim Center for Finance and the
Department of Operations Research \& Financial Engineering at
Princeton University. Grateful acknowledgement is made for
hospitality.}


\end{document}